\documentclass[11pt, a4paper]{article}
\usepackage{}
\usepackage{amsfonts}
\usepackage{mathbbold}
\usepackage{bbm}
\usepackage{mathrsfs}
\usepackage{latexsym}
\usepackage{graphicx}
\usepackage{amssymb}

\newtheorem{theorem}{Theorem}[section]
\newtheorem{lemma}[theorem]{Lemma}

\title{{\Large \bf Further results on the least $Q$-eigenvalue of a graph with fixed domination number\thanks{Supported by NSFC
(Nos. 11771376 \& 11571252), ``333" Project of  Jiangsu (2016) \& KPPAT of Anhui (JXBJZD 2016082).}}}
\author{Guanglong Yu$^a$\thanks{E-mail addresses:
yglong01@163.com.}
 ~ Yarong Wu$^b$ ~ Mingqing Zhai$^c$ ~
\\ ~ \\
{\footnotesize $^a$Department of Mathematics, Lingnan normal
nniversity,  Zhanjiang, 524048, Guangdong, China}\\
{\footnotesize $^b$SMU college of art and science, Shanghai maritime
university, Shanghai, 200135, China}\\
{\footnotesize $^c$ School of mathematics and finance, Chuzhou university, Chuzhou, 239000, Anhui, China}}
\date{}

\begin{document}
\maketitle

\begin{abstract}
In this paper, we proceed on determining the minimum $q_{min}$ among the connected nonbipartite graphs on $n\geq 5$ vertices and with domination number $\frac{n+1}{3}<\gamma\leq \frac{n-1}{2}$. Further results obtained are as follows:

$\mathrm{(i)}$ among all nonbipartite connected graph of order $n\geq 5$ and with domination number $\frac{n-1}{2}$, the minimum $q_{min}$ is completely determined;

$\mathrm{(ii)}$ among all nonbipartite graphs of order $n\geq 5$, with odd-girth $g_{o}\leq5$ and domination number at least $\frac{n+1}{3}<\gamma\leq \frac{n-2}{2}$, the minimum $q_{min}$ is completely determined.


\end{abstract}
\baselineskip 18.6pt

\section{Introduction}

\ \ \ \ All graphs considered in this paper are connected, undirected and
simple, i.e., no loops or multiple edges are allowed.
We denote by $\parallel S\parallel$ the $cardinality$ of a set $S$,
and denote by $G=G[V(G)$, $E(G)]$ a graph with vertex set
$V(G)=\{v_1, v_2, \ldots, v_n\}$ and edge set $E(G)$ where $\parallel V(G)\parallel= n$ is the $order$
and $\parallel E(G)\parallel= m$ is the $size$.

In a graph, if vertices $v_{i}$ and $v_{j}$ are adjacent (denoted by $v_{i}\sim v_{j}$), we say that they $dominate$ each other.
A vertex set $D$ of a graph $G$ is said to be a $dominating$ $set$ if every vertex of $V(G)\setminus D$ is
adjacent to (dominated by) at least one vertex in $D$. The $domination$ $number$ $\gamma(G)$ ($\gamma$, for short) is the
minimum cardinality of all dominating sets of $G$. For a graph $G$, a dominating set is called a $minimal$ $dominating$ $set$ if its cardinality is $\gamma(G)$. A well known result about $\gamma(G)$ is that for a graph $G$
of order $n$ containing no isolated vertex, $\gamma\leq \frac{n}{2}$ \cite{ORE}. A comprehensive
study of issues relevant to dominating set of a graph has been undertaken because of its good applications \cite{EWB}, \cite{YWWY}.

Recall
that $Q(G)= D(G) + A(G)$ is called the $signless$ $Laplacian$ $matrix$ (or $Q$-$matrix$) of $G$, where $D(G)=
\mathrm{diag}(d_{1}, d_{2},
\ldots, d_{n})$ with $d_{i}= d_{eg}(v_{i})$ being the degree of
vertex $v_{i}$ $(1\leq i\leq n)$, and $A(G)$ is the adjacency matrix of $G$. The signless
Laplacian has attracted the
attention of many researchers and it is
being promoted by many researchers \cite{CCRS}, \cite{D.P.S}-\cite{LOA}, \cite{WF}.

The least eigenvalue of $Q(G)$,
denote by $q_{min}(G)$ or $q_{min}$, is called the $least$ $Q$-$eigenvalue$ of $G$. Because $Q(G)$ is positive semi-definite,
we have $q_{min}(G)\geq 0.$
From \cite{D.P.S}, we know that, for a connected
graph $G$, $q_{min}(G)= 0$ if and only if $G$ is bipartite.
Consequently, in \cite{DR}, $q_{min}$ was studied as a measure of nonbipartiteness of a graph. One can notice
that there are quite a few results about $q_{min}$. In \cite{CCRS}, D.M. Cardoso et al. determined the graphs
with the the minimum $q_{min}$ among all the connected nonbipartite
graphs with a prescribed number of vertices. In \cite{LOA}, L. de Lima et al. surveyed some known results about $q_{min}$ and also presented some
new results. In \cite{FF}, S. Fallat, Y. Fan investigated the relations
between $q_{min}$ and some parameters reflecting
the graph bipartiteness. In \cite{WF}, Y. Wang, Y. Fan investigated $q_{min}$
of a graph under some perturbations, and minimized $q_{min}$ among the connected graphs with fixed order
which contains a given nonbipartite graph as an induced subgraph. Recently, in \cite{RZSG}, the authors determined all non-bipartite hamiltonian graphs whose $q_{min}$ attains the minimum.

Recall that a $lollipop$ $graph$ $L_{g, l}$ is a graph composed of a cycle $\mathbb{C}=v_{1}v_{2}\cdots v_{g}v_{1}$ and a path $\mathbb{P}=v_{g}v_{g+1}\cdots v_{g+l}$ with
$l\geq 1$. For given $g$ and $l$, a graph of order $n$ is called a $F_{g, l}$-$graph$ if it
is obtained by attaching $n-g-l$ pendant vertices to some nonpendant vertices of a $L_{g, l}$. If $l=1$, a $F_{g, l}$-graph is also called a $sunlike$ graph. In a graph, a vertex is called a $p$-$dominator$ (or $support$ $vertex$) if it dominates a pendant vertex.
In a $F_{g, l}$-graph if each $p$-dominator other than $v_{g+l-1}$ is  attached with exactly one pendant vertex, then this graph is called a $\mathcal {F}_{g, l}$-$graph$. A $\mathcal {F}_{g, l}$-graph is called a $\mathcal {F}^{\circ}_{g, l}$-$graph$ if $v_{g}$ is a $p$-dominator. In the following paper, for unity, for a $\mathcal {F}_{g, l}$-graph, $\mathbb{C}$ and $\mathbb{P}$ are expressed as above.

\setlength{\unitlength}{0.6pt}
\begin{center}
\begin{picture}(629,207)
\put(1,190){\circle*{4}}
\put(1,139){\circle*{4}}
\qbezier(1,190)(1,165)(1,139)
\put(31,166){\circle*{4}}
\qbezier(1,190)(16,178)(31,166)
\qbezier(1,139)(16,153)(31,166)
\put(85,166){\circle*{4}}
\put(53,166){\circle*{4}}
\put(368,167){\circle*{4}}
\put(110,166){\circle*{4}}
\put(101,166){\circle*{4}}
\put(93,166){\circle*{4}}
\put(189,166){\circle*{4}}
\put(254,167){\circle*{4}}
\put(360,167){\circle*{4}}
\qbezier(254,167)(307,167)(360,167)
\qbezier(110,166)(149,166)(189,166)
\put(135,166){\circle*{4}}
\put(157,166){\circle*{4}}
\put(157,189){\circle*{4}}
\qbezier(157,166)(157,178)(157,189)
\put(218,166){\circle*{4}}
\qbezier(189,166)(203,166)(218,166)
\put(229,167){\circle*{4}}
\put(238,167){\circle*{4}}
\put(247,167){\circle*{4}}
\put(334,167){\circle*{4}}
\put(306,166){\circle*{4}}
\put(280,167){\circle*{4}}
\put(306,190){\circle*{4}}
\qbezier(306,166)(306,178)(306,190)
\put(378,167){\circle*{4}}
\put(387,167){\circle*{4}}
\put(471,167){\circle*{4}}
\put(444,167){\circle*{4}}
\put(421,167){\circle*{4}}
\put(444,191){\circle*{4}}
\put(395,167){\circle*{4}}
\put(497,167){\circle*{4}}
\qbezier(395,167)(446,167)(497,167)
\qbezier(444,167)(444,179)(444,191)
\put(76,166){\circle*{4}}
\put(592,167){\circle*{4}}
\put(560,167){\circle*{4}}
\put(629,167){\circle*{4}}
\qbezier(560,167)(594,167)(629,167)
\put(507,167){\circle*{4}}
\put(516,167){\circle*{4}}
\put(525,167){\circle*{4}}
\put(533,167){\circle*{4}}
\put(541,167){\circle*{4}}
\put(549,167){\circle*{4}}
\qbezier(31,166)(53,166)(76,166)
\put(29,86){\circle*{4}}
\put(29,43){\circle*{4}}
\put(68,64){\circle*{4}}
\put(118,64){\circle*{4}}
\put(163,64){\circle*{4}}
\put(289,64){\circle*{4}}
\put(574,84){\circle*{4}}
\put(578,78){\circle*{4}}
\put(581,72){\circle*{4}}
\put(234,64){\circle*{4}}
\put(412,64){\circle*{4}}
\put(176,64){\circle*{4}}
\put(185,64){\circle*{4}}
\put(194,64){\circle*{4}}
\put(204,64){\circle*{4}}
\put(214,64){\circle*{4}}
\put(223,64){\circle*{4}}
\qbezier(29,86)(29,65)(29,43)
\qbezier(29,86)(48,75)(68,64)
\qbezier(29,43)(48,54)(68,64)
\qbezier(234,64)(323,64)(412,64)
\qbezier(68,64)(115,64)(163,64)
\put(-9,196){$v_{2}$}
\put(27,172){$v_{3}$}
\put(45,155){$v_{4}$}
\put(-15,133){$v_{1}$}
\put(151,154){$v_{a_{1}}$}
\put(297,154){$v_{a_{2}}$}
\put(436,154){$v_{a_{k}}$}
\put(335,51){$v_{\varepsilon-k-2}$}
\put(271,50){$v_{\varepsilon-k-3}$}
\put(24,94){$v_{2}$}
\put(65,72){$v_{3}$}
\put(23,32){$v_{1}$}
\put(226,-15){Fig. 1.1. $\mathcal {H}^{k}_{1}$, $\mathcal {H}^{k}_{2}$}
\put(293,115){$\mathcal {H}^{k}_{1}$}
\put(297,14){$\mathcal {H}^{k}_{2}$}
\put(626,153){$v_{\varepsilon}$}
\put(579,154){$v_{\varepsilon-1}$}
\put(348,64){\circle*{4}}
\put(534,64){\circle*{4}}
\put(477,64){\circle*{4}}
\put(586,64){\circle*{4}}
\put(424,64){\circle*{4}}
\put(433,64){\circle*{4}}
\put(442,64){\circle*{4}}
\put(450,64){\circle*{4}}
\put(458,64){\circle*{4}}
\put(466,64){\circle*{4}}
\qbezier(477,64)(531,64)(586,64)
\put(526,51){$v_{\varepsilon-1}$}
\put(585,52){$v_{\varepsilon}$}
\put(611,194){\circle*{4}}
\qbezier(592,167)(601,181)(611,194)
\put(624,174){\circle*{4}}
\put(619,180){\circle*{4}}
\put(615,186){\circle*{4}}
\qbezier(1,139)(1,165)(1,190)
\qbezier(1,139)(1,165)(1,190)
\qbezier(1,139)(1,165)(1,190)
\qbezier(1,139)(1,165)(1,190)
\put(571,92){\circle*{4}}
\qbezier(534,64)(552,78)(571,92)
\put(477,85){\circle*{4}}
\qbezier(477,64)(477,75)(477,85)
\put(412,85){\circle*{4}}
\qbezier(412,64)(412,75)(412,85)
\put(348,85){\circle*{4}}
\qbezier(348,64)(348,75)(348,85)
\put(397,49){$v_{\varepsilon-k-1}$}
\put(466,49){$v_{\varepsilon-2}$}
\end{picture}
\end{center}

Let $\mathcal {H}^{k}_{1}$ be a $\mathcal {F}_{3, \varepsilon-3}$-graph of order $n\geq 4$ where there are $k\geq 0$ $p$-dominators among $v_{1}$, $v_{2}$, $\ldots$, $\varepsilon-2$ ($\varepsilon\geq3$. see Fig. 1.1). If $k\geq 1$, in $\mathcal {H}^{k}_{1}$, suppose $v_{a_{j}}s$ are $p$-dominators where $1\leq j\leq k$, $1\leq a_{1}<a_{2}<\cdots <a_{k}\leq \varepsilon-2$, and suppose $v_{\tau_{j}}$ is the pendant vertex attached to $v_{a_{j}}$. Let $\mathcal {H}^{k}_{2}=\mathcal {H}^{k}_{1}-\sum\limits_{j=1}^{k} v_{\tau_{j}}v_{a_{j}}+\sum\limits_{j=1}^{k} v_{\tau_{j}}v_{\varepsilon-2-k+j}$ (see Fig. 1.1). If $k= 0$, then $\mathcal {H}^{0}_{1}=\mathcal {H}^{0}_{2}$. If $\alpha\geq 1$, we denoted by $\mathscr{H}_{3,\alpha}$ the graph $\mathcal {H}^{\alpha-1}_{2}$ of order $n$ in which there are $\alpha$ $p$-dominators and $v_{\varepsilon-1}$ has only one pendant vertex (where $\varepsilon=n-\alpha+1$); if $\alpha= 0$, we let $\mathscr{H}_{3,0}=C_{3}=v_{1}v_{2}v_{3}v_{1}$.

In \cite{YFT} and \cite{YGZW}, the authors first considered the relation between $q_{min}$ of a graph and its domination number. Among all the nonbipartite graphs with both order $n\geq 4$ and domination number $\gamma\leq \frac{n+1}{3}$, they characterized the graphs with the minimum $q_{min}$. A remaining open problem is that how about the $q_{min}$ of the connected nonbipartite graph on $n$ vertices with domination number $\frac{n+1}{3}<\gamma\leq \frac{n}{2}$. In \cite{YMCG}, the authors proceeded on considering this problem. Among the nonbipartite graphs of order $n=4$, the minimum $q_{min}$ is completely determined; among the nonbipartite graphs of order $n$ and with given domination number $\frac{n}{2}$, the minimum $q_{min}$ is completely determined; further results about the domination number, the $q_{min}$ of a graph as well as their relation are represented. An open problem still left is that how to determine the minimum $q_{min}$ of the connected nonbipartite graph on $n\geq 5$ vertices with domination number $\frac{n+1}{3}<\gamma\leq \frac{n-1}{2}$. Let $\mathbb{S}=\mathscr{H}_{3,\alpha}$ be of order $n\geq 4$ where $\alpha$ is the least integer such that $\lceil\frac{n-2\alpha-2}{3}\rceil+\alpha=\gamma$. In \cite{YMCG}, the authors represented some structural characterizations about the minimum $q_{min}$ for this problem, and conjectured that such $\mathbb{S}$ has the smallest $q_{min}$. However, the problem seems really difficult to solve. Motivated by proceeding on solving this problem, we go on with our research and get some further results as follows.

\begin{theorem}\label{th05,3} 
Let $G$ be a nonbipartite connected graph of order $n\geq 5$ and with domination number $\frac{n-1}{2}$. Then $q_{min}(G)\geq q_{min}(\mathscr{H}_{3,\frac{n-1}{2}})$ with equality if and only if $G\cong \mathscr{H}_{3,\frac{n-1}{2}}$.
\end{theorem}

\begin{theorem}\label{th05,4} 
Among all nonbipartite graphs of order $n\geq 5$, with odd-girth $g_{o}\leq5$ (length of the shortest odd cycle in this graph) and domination number $\frac{n+1}{3}<\gamma\leq \frac{n-2}{2}$, then the least $q_{min}$ attains the minimum uniquely at a $\mathscr{H}_{3,\alpha}$ where $\alpha\leq\frac{n-3}{2}$ is the least integer such that $\lceil\frac{n-2\alpha-2}{3}\rceil+\alpha=\gamma$.

\end{theorem}

\section{Preliminary}

\ \ \ \ In this section, we introduce some notations and some working lemmas.

Denote by $P_n$, $C_{n}$, $K_{n}$, a $path$, a $n$-cycle (of length $n$), a $complete$ graph of order $n$ respectively. If $k$ is odd, we say $C_{k}$ an $odd$ $cycle$. The $girth$ of a graph $G$, denoted by $g$, is the length of the shortest cycle in $G$. The $odd$-$girth$ for a nonbipartite graph $G$,
denoted by $g_{o}(G)$ or $g_{o}$, is the length of the shortest odd cycle in this graph.  $G-v_{i}v_{j}$
 denotes the graph obtained from $G$ by deleting the edge $v_{i}v_{j}\in
 E(G)$, and let $G-v_{i}$
 denote the graph obtained from $G$ by deleting the vertex $v_{i}$ and the edges incident with $v_{i}$.
 Similarly,  $G+v_{i}v_{j}$ is the graph obtained from $G$ by adding an edge $v_{i}v_{j}$ between its two nonadjacent vertices $v_{i}$
 and $v_{j}$. Given an vertex set $S$, $G-S$ denotes the graph obtained by deleting all the vertices in $S$ from $G$ and the edges incident with any vertex in $S$.

A connected graph $G$ of order $n$ is called a $unicyclic$ graph if
$\|E(G)\|=n$. For $S\subseteq V(G)$, let $G[S]$
denote the subgraph induced by $S$.
Denoted by $d_{istG}(v_{i}, v_{j})$ the $distance$ between two vertices
$v_{i}$ and $v_{j}$ in a graph $G$.

For a graph $G$ of order $n$, let $X=(x_1, x_2, \ldots, x_n)^T \in R^n$ be defined on $V(G)$, i.e.,
each vertex $v_i$ is mapped to the entry $x_i$; let $|x_i|$ denote the $absolute$ $value$ of $x_i$.
One can find
that $X^TQ(G)X =\sum_{v_iv_j\in E(G)}(x_{i} + x_{j})^2.$
In addition, for an arbitrary unit vector $X\in R^n$, $q_{min}(G) \le X^TQ(G)X$,
with equality if and only if $X$ is an eigenvector corresponding to $q_{min}(G)$.

\begin{lemma}{\bf \cite{D.R.S}} \label{le02,01} 
Let $G$ be a graph on $n$ vertices and $m$ edges, and let $e$ be an
edge of $G$. Let $q_{1}\geq q_{2}\geq \cdots \geq q_{n}$ and
$s_{1}\geq s_{2}\geq \cdots \geq s_{n}$ be the $Q$-eigenvalues of  $G$
and $G-e$ respectively. Then $0\leq s_{n}\leq q_{n}\leq \cdots \leq
s_{2}\leq q_{2}\leq s_{1}\leq q_{1}.$
\end{lemma}

Let $G_1$ and $G_2$ be two disjoint graphs, and let $v_1\in V(G_1),$  $v_2\in V(G_2)$. The $coalescence$ of $G_1$ and $G_2$,
denoted by $G_1(v_1)\diamond G_2(v_2)$ or $G_1(u)\diamond G_2(u)$, is obtained from $G_1$, $G_2$ by identifying $v_1$ with $v_2$ and forming a new vertex
$u$ where for $i=1, 2$, $G_{i}$ can be trivial (that is, $G_{i}$ is only one vertex). For a connected
graph $G = G_1(u)\diamond G_2(u)$, $i=1$, $2$, $G_i$ is
called a
$branch$ of $G$ with root $u$. For a vector $X=(x_1, x_2, \ldots, x_n)^T \in R^n$ defined on $V(G)$, a branch $H$ of $G$ is called a $zero$ $branch$ with respect to $X$ if $x_{i} = 0$  for all $v_{i} \in V(H)$; otherwise, it
is called a $nonzero$ $branch$ with respect to $X$.

\begin{lemma}{\bf \cite{WF}}\label{le02,02} 
Let $G$ be a connected graph which contains a bipartite branch $H$ with root $v_{s}$, and let $X$ be an eigenvector of $G$
corresponding
to $q_{min}(G)$.

{\normalfont (i)} If $x_{s} = 0$, then $H$ is a zero branch of G with respect to $X$;

{\normalfont (ii)} If $x_{s}\neq 0$, then $x_{p}\neq 0$ for every vertex $v_{p}\in V(H)$. Furthermore, for every vertex $v_{p}\in V(H)$,
$x_{p}x_{s}$ is either positive or negative depending on whether $v_{p}$ is or is not in the same part of the bipartite graph $H$ as
$v_{s}$; consequently, $x_{p}x_{t} < 0$ for each edge $v_{p}v_{t} \in E(H)$.
\end{lemma}

\begin{lemma}{\bf \cite{WF}}\label{le02,03} 
Let $G$ be a connected nonbipartite graph of order $n$, and let $X$ be an eigenvector of $G$ corresponding to $q_{min}(G)$.
$T$ is a tree which is a nonzero branch of $G$ with respect to $X$ and with root $v_{s}$. Then $|x_{t}| < |x_{p}|$
 whenever $v_{p},$ $v_{t}$
are vertices of $T$ such that $v_{t}$ lies on the unique path from $v_{s}$ to $v_{p}$.
\end{lemma}

\begin{lemma}{\bf \cite{YGX}}\label{le0,04} 
 Let $G = G_1(v_2) \diamond T(u)$ and $G^* = G_1(v_1)\diamond T(u)$, where $G_1$ is a connected nonbipartite
graph containing two distinct vertices $v_1, v_2$, and $T$ is a nontrivial tree. If there exists an
 eigenvector $X=(\,x_1$, $x_2$, $\ldots$, $x_k$, $\ldots)^T$ of $G$ corresponding to $q_{min}(G)$ such that
 $|x_1| > |x_2|$ or $|x_1| = |x_2| > 0$, then $q_{min}(G^*)<q_{min}(G)$.
\end{lemma}

\begin{lemma}{\bf \cite{YGX}}\label{le0,05} 
Let $G = C(v_0)\diamond B(v_0)$ be a graph of order $n$, where $C=v_0$$v_1$$v_2$ $\cdots$ $v_{2k}$ is a cycle of length $2k+1$, and $B$ is a bipartite graph of order $n-2k$. Then there exists an  eigenvector
$X=(\,x_{0}$, $x_{1}$, $x_{2}$, $\ldots$, $x_{2k}$ $\,)^T$
corresponding to $q_{min}(G)$ satisfying the following:

{\normalfont (i)} $|x_{0}|=\max\{|x_{i}|\,|\, v_{i}\in V(C)\}>0$;

{\normalfont (ii)}  $x_{i}=x_{2k-i+1}$ for $i=1, 2, \ldots, k$;

{\normalfont (iii)} $x_{i}x_{i-1}\leq 0$ for $i=1, 2, \ldots, k$, $x_{2k}x_{0}\leq 0$ and $x_{2k-i+1}x_{2k-i+2}\leq 0$ for $i=2, \ldots, k$.

Moreover, if $2k+1<n$, then the multiplicity of $q_{min}(G)$ is one, and then any  eigenvector
corresponding to $q_{min}(G)$ satisfies {\normalfont (i), (ii), (iii)}.
\end{lemma}

\begin{lemma}{\bf \cite{KDAS}}\label{le0,06} 
Let $G$ be a connected graph of order $n$. Then
$q_{min}<\delta$, where $\delta$ is the minimal vertex degree of $G$.
\end{lemma}

\begin{lemma}{\bf \cite{YGZW}}\label{le0,07} 
Let $G$ be a nonbipartite graph with domination number $\gamma(G)$. Then $G$ contains a nonbipartite
unicyclic spanning subgraph $H$
with both $g_{o}(H)=g_{o}(G)$ and $\gamma(H)=\gamma(G)$.
\end{lemma}

\begin{lemma}{\bf \cite{YGZW}}\label{le0,08} 
Suppose a graph $G$ contains pendant vertices. Then

$\mathrm{(i)}$ there must be a minimal dominating set of $G$ containing
all of its $p$-dominators but no any pendant vertex;

$\mathrm{(ii)}$ if $v$ is a $p$-dominator of $G$ and at least two pendant vertices are adjacent to $v$, then any
minimal dominating set of $G$ contains $v$ but no any pendant vertex adjacent to $v$.
\end{lemma}

\begin{lemma}{\bf \cite{MHM}}\label{le0,09} 
(i) For a path $P_{n}$, we have
$\gamma(P_{n})=\lceil\frac{n}{3}\rceil$.

(ii) For a cycle $C_{n}$, we have
$\gamma(C_{n})=\lceil\frac{n}{3}\rceil$.
\end{lemma}

We define the corona $G$ of graphs $G_{1}$ and $G_{2}$ as follows. The $corona$ $G = G_{1}\circ G_{2}$ is the graph
formed from one copy of $G_{1}$ and $\parallel V (G_{1})\parallel$ copies of $G_{2}$ where the $i$th vertex of $G_{1}$ is adjacent to
every vertex in the $i$th copy of $G_{2}$.

\begin{lemma}{\bf \cite{CPNH}}\label{le0,10} 
Let $G$ be a graph of order $n$. $\gamma(G) =\frac{n}{2}$ if
and only if the components of $G$ are the cycle $C_{4}$ or the corona $H\circ K_{1}$ for any connected graph $H$.
\end{lemma}

Denote by $C_{3,\, k}^*$ the graph obtained by attaching a $C_{3}$ to an end vertex of a path of length $k$ and attaching $n-3-k$ pendant vertices to the other end vertex of this path.

\begin{lemma}{\bf \cite{YGZW}}\label{th2,18} 
Among all the nonbipartite graphs with both order $n\geq 4$ and domination number $\gamma\leq \frac{n+1}{3}$, we have

\normalfont (i) if $n=3\gamma-1$, $3\gamma$, $3\gamma+1$, then the graph with the minimal least $Q$-eigenvalue attains uniquely
at  $C_{3,\, n-4}^*$;

\normalfont (ii) if $n\geq3\gamma+2$, then the graph with the minimal least $Q$-eigenvalue attains uniquely at
$C_{3,\, 3\gamma-3}^*$.
\end{lemma}

\begin{lemma}{\bf \cite{YMCG}}\label{th03,02} 
Among all nonbipartite unicyclic graphs of order $n$, and with both domination number $\gamma$ and girth $g$ ($g\leq n-1$), the
minimum $q_{min}$ attains at a $\mathcal {F}_{g, l}$-graph $G$ for some $l$. Moreover, for this graph $G$, suppose that $X=(x_1$, $x_2$, $x_3$, $\ldots$, $x_n)^T$ is a unit
eigenvector corresponding to $q_{min}(G)$. Then we have that $| x_{g}|> 0$, and $|x_{g+l-1}|=\max\{|x_{i}|\mid v_{i}$ is a $p$-dominator$\,\}$.
\end{lemma}

 In $\mathcal {H}^{k}_{2}$, for $j=1$, $2$, $\ldots$, $k$, suppose $v_{\tau_{\varepsilon-2-k+j}}$ is the pendant vertex attached to vertex $v_{\varepsilon-2-k+j}$. Suppose $v_{\omega_{1}}$, $v_{\omega_{2}}$, $\ldots$, $v_{\omega_{s}}$ are the pendant vertices attached to vertex $v_{\varepsilon-1}$. If $s\geq 2$, let $\mathcal {H}^{k}_{3}=\mathcal {H}_{2}^{k}-v_{\varepsilon-1-k}v_{\tau_{\varepsilon-1-k}}+
v_{\varepsilon-1}v_{\tau_{\varepsilon-1-k}}-\sum\limits_{j=2}^{s}v_{\varepsilon-1}v_{\omega_{j}}+
\sum\limits_{j=2}^{s}v_{\omega_{1}}v_{\omega_{j}}$.  Let
$\mathcal {H}^{k-1}_{4}=\mathcal {H}^{k}_{2}-v_{\varepsilon-1-k}v_{\tau_{\varepsilon-1-k}}+v_{\varepsilon-1}v_{\tau_{\varepsilon-1-k}}$, $\mathcal {H}^{k-2}_{5}=\mathcal {H}^{k-1}_{4}-v_{\varepsilon-k}v_{\tau_{\varepsilon-k}}+v_{\varepsilon-1}v_{\tau_{\varepsilon-k}}$.

\begin{lemma}{\bf \cite{YMCG}}\label{th03,05} 
\

$\mathrm{(i)}$ $\gamma(\mathcal {H}^{k}_{1})\leq\gamma(\mathcal {H}^{k}_{2})$.

$\mathrm{(ii)}$ If $\varepsilon-k-1\leq 2$, then $\gamma(\mathcal {H}^{k}_{2})=k+1$ and $\gamma(\mathcal {H}^{k-1}_{4})= \gamma(\mathcal {H}^{k}_{2})-1$;

$\mathrm{(iii)}$ If $\varepsilon-k-1\geq 3$, then $\gamma(\mathcal {H}^{k}_{2})=\lceil\frac{\varepsilon-k-4}{3}\rceil+k+1$;

$\mathrm{(iv)}$ $\gamma(\mathcal {H}^{k}_{2})\leq \gamma(\mathcal {H}^{k}_{3})$;

$\mathrm{(v)}$ If $\varepsilon-k-1\geq 3$, $\frac{\varepsilon-k-4}{3}\neq t$ where $t$ is a nonnegative integral number, then $\gamma(\mathcal {H}^{k-1}_{4})= \gamma(\mathcal {H}^{k}_{2})-1$;

$\mathrm{(vi)}$ If $\varepsilon-k-1\geq 3$, $\frac{\varepsilon-k-4}{3}= t$ where $t$ is a nonnegative integral number, $\gamma(\mathcal {H}^{k-1}_{4})= \gamma(\mathcal {H}^{k}_{2})$, $\gamma(\mathcal {H}^{k-2}_{5})= \gamma(\mathcal {H}^{k}_{2})-1$.
\end{lemma}

\begin{lemma}{\bf \cite{YMCG}}\label{cl03,06} 
\

$\mathrm{(i)}$ $\gamma(\mathscr{H}_{3,0})=1$;

$\mathrm{(ii)}$ If $\alpha\geq 1$ and $n-2\alpha\leq 2$, then $\gamma(\mathscr{H}_{3,\alpha})= \alpha$;

$\mathrm{(iii)}$ If $\alpha\geq 1$ and $n-2\alpha\geq 3$, then $\gamma(\mathscr{H}_{3,\alpha})=\lceil\frac{n-2\alpha-2}{3}\rceil+\alpha$.
\end{lemma}

\section{Domination number and the structure of a graph}

Let $G^{\ast}$ be a sunlike graph of order $n$ and with both girth $g$ and $k$ $p$-dominators $v_{1}$, $v_{2}$, $\ldots$, $v_{k}$ on $\mathbb{C}$.
\begin{lemma}\label{le03,07} 
Let $G$ be a sunlike graph of order $n$ and with both girth $g$ and $k$ $p$-dominators on $\mathbb{C}$. Then $\gamma(G)\leq \gamma(G^{\ast})$, where $\gamma(G^{\ast})=k+\lceil\frac{g-k-2}{3}\rceil$.
\end{lemma}

\begin{proof}
Suppose $v_{i_{1}}$, $v_{i_{2}}$, $\ldots$, $v_{i_{k}}$ are the $k$ $p$-dominators on $\mathbb{C}$ in $G$, where $1\leq i_{1}<i_{2}<\cdots< i_{k}\leq g$. Suppose that there exists some $1\leq z\leq k$ such that $i_{z+1}-i_{z}\geq 2$, where if $z=k$, we let $i_{k+1}=i_{1}$ and $i_{k+1}-i_{k}=i_{1}+g-i_{k}$. Let $H=G-\sum^{i_{z+1}-1}_{s=i_{z}+1}v_{s}$.

{\bf Assertion 1} If $i_{z+1}-i_{z}\leq 3$, then $\gamma(H)=\gamma(G)$. By Lemma \ref{le0,08}, there is a minimal dominating set $D$ of $G$ which contains all the $k$ $p$-dominators but no any pendant vertex. Thus both $v_{i_{z+1}}$ and $v_{i_{z}}$ are in $D$. Note the minimality of $D$ and $2\leq i_{z+1}-i_{z}\leq 3$. Then $D\cap \{v_{i_{z}+1}\}=\emptyset$ if $i_{z+1}-i_{z}=2$; $D\cap \{v_{i_{z}+1},v_{i_{z+1}-1}\}=\emptyset$ if $i_{z+1}-i_{z}=3$. Thus $D$ is also a dominating set of $H$. This implies that $\gamma(H)\leq \gamma(G)$. Note that for $H$, by Lemma \ref{le0,08}, there is a minimal dominating set $D'$ which contains all the $k$ $p$-dominators but no any pendant vertex. Thus both $v_{i_{z+1}}$ and $v_{i_{z}}$ are in $D'$. Then $v_{i_{z}+1}$ is dominated by $D'$ if $i_{z+1}-i_{z}=2$; $v_{i_{z}+1}$, $v_{i_{z+1}-1}$ are is dominated by $D'$ if $i_{z+1}-i_{z}=3$. Consequently, $D'$ is also a dominating set of $G$. This implies that $\gamma(G)\leq\gamma(H)$. As a result, it follows that $\gamma(H)=\gamma(G)$. And then our assertion holds.

{\bf Assertion 2} If $i_{z+1}-i_{z}\geq 4$, then $\gamma(G)=\gamma(H)+\gamma(P_{i_{z},i_{z+1}})$ where $P_{i_{z},i_{z+1}}=v_{i_{z}+2}v_{i_{z}+3}\cdots$ $v_{i_{z+1}-2}$. By Lemma \ref{le0,08}, there is a minimal dominating set $D$ of $G$ which contains all the $k$ $p$-dominators but no any pendant vertex. Thus both $v_{i_{z+1}}$ and $v_{i_{z}}$ are in $D$. We claim that at most one of $v_{i_{z}+1}$, $v_{i_{z}+2}$ is in $D$. Otherwise, suppose that both $v_{i_{z}+1}$ and $v_{i_{z}+2}$ are in $D$. Then $D\setminus \{v_{i_{z}+1}\}$ is also a dominating set of $G$, which contradicts the minimality of $D$. Consequently, our claim holds. Similarly, we get that at most one of $v_{i_{z+1}-2}$, $v_{i_{z+1}-1}$ is in $D$. Thus we let $D^{\circ}=((D\cup\{v_{i_{z}+2}, v_{i_{z+1}-2}\})\setminus \{v_{i_{z}+1}, v_{i_{z+1}-1}\})\cap V(P_{i_{z},i_{z+1}})$ if $v_{i_{z}+1}\in D$, $v_{i_{z+1}-1}\in D$;
let $D^{\circ}=((D\cup\{v_{i_{z}+2}\})\setminus \{v_{i_{z}+1}\})\cap V(P_{i_{z},i_{z+1}})$ if $v_{i_{z}+1}\in D$ and $v_{i_{z+1}-1}\notin D$; let $D^{\circ}=((D\cup\{v_{i_{z+1}-2}\})\setminus \{v_{i_{z+1}-1}\})\cap V(P_{i_{z},i_{z+1}})$ if $v_{i_{z}+1}\notin D$ and $v_{i_{z+1}-1}\in D$; let $D^{\circ}=(D\cap V(P_{i_{z},i_{z+1}})$ if $v_{i_{z}+1}\notin D$ and $v_{i_{z+1}-1}\notin D$. Note that $D^{\ast}=D\setminus (V(P_{i_{z},i_{z+1}})\cup \{v_{i_{z}+1}$, $v_{i_{z+1}-1}\})$ is a dominating set of $H$, $D^{\circ}\cup D^{\ast}$ is a dominating set of $G$ with cardinality $\gamma(G)$, and note that $D^{\circ}$ is a dominating set of $P_{i_{z},i_{z+1}}$. Thus $\gamma(P_{i_{z},i_{z+1}})\leq\parallel D^{\circ}\parallel$. Note that both $v_{i_{z+1}-1}$ and $v_{i_{z}+1}$ are dominated by $D^{\ast}$. Consequently, for any minimal dominating set $B$ of $P_{i_{z},i_{z+1}}$, then $B\cup D^{\ast}$ is also a dominating set of $G$. Note that $\parallel B\parallel =\gamma(P_{i_{z},i_{z+1}})\leq\parallel D^{\circ}\parallel$.
As a result, $\parallel B\cup D^{\ast}\parallel\leq \parallel D\parallel=\gamma(G)$. Note that the minimality of $D$. Then $\parallel D^{\circ}\parallel=\parallel B\parallel =\gamma(P_{i_{z},i_{z+1}})$, and then it follows that $\gamma(G)=\gamma(H)+\gamma(P_{i_{z},i_{z+1}})$.

Denote by $\tau_{i_{j},i_{j+1}}$ the dominating index where we let $i_{k+1}=i_{1}$ if $i=k$. Let $\tau_{i_{j},i_{j+1}}=0$ if $i_{j+1}-i_{j}\leq 3$; let $\tau_{i_{j},i_{j+1}}=\gamma(P_{i_{j},i_{j+1}})$ if $i_{j+1}-i_{j}\geq 4$. Thus from Assertion 1, Assertion 2 and Lemma \ref{le0,08}, we get that $\gamma(G)=k+\sum^{k}_{i=1}\tau_{i_{j},i_{j+1}}$. By Lemma \ref{le0,09}, it follows that $\tau_{i_{j},i_{j+1}}=\gamma(P_{i_{j},i_{j+1}})=\lceil\frac{i_{j+1}-i_{j}-3}{3}\rceil$ if $i_{j+1}-i_{j}\geq 4$. Note that for any two nonnegative integers $x$ and $y$, we have $\lceil\frac{x}{3}\rceil+\lceil\frac{y}{3}\rceil\leq \lceil\frac{x+y}{3}\rceil$. Then $$\displaystyle\sum^{k}_{i=1}\tau_{i_{j},i_{j+1}}=\sum_{\tau_{i_{s},i_{s+1}}\neq 0}\tau_{i_{s},i_{s+1}}\leq \left\lceil\frac{\sum_{\tau_{i_{s},i_{s+1}}\neq 0}(i_{s+1}-i_{s}-3)}{3}\right\rceil\leq \left\lceil\frac{g-k-2}{3}\right\rceil.$$
Thus $\gamma(G)\leq k+\lceil\frac{g-k-2}{3}\rceil$. Noting that by Assertion1 and Assertion 2, we have $\gamma(G^{\ast})=k+\lceil\frac{g-k-2}{3}\rceil$. Then the result follows as desired. This completes the proof. \ \ \ \ \ $\Box$
\end{proof}

\begin{theorem}\label{le03,08} 
Suppose that $\mathcal {G}$ is a nonbipartite $\mathcal {F}_{g, l}$-graph with $\gamma(\mathcal {G}) =\frac{n-1}{2}$, $g\geq 5$ and order $n\geq g+1$, and suppose there are exactly $f$ vertices of the unique cycle $\mathbb{C}$ such that none of them is $p$-dominator. Then we get

$\mathrm{(i)}$ if $f=g$, then $g=5$;

$\mathrm{(ii)}$ if $f\neq g$, then $f\leq 3$ and $f\neq2$;

$\mathrm{(iii)}$ if $f= 3$, then the three vertices are consecutive on $\mathbb{C}$, i.e., they are $v_{i-1}$, $v_{i}$, $v_{i+1}$ for some $1\leq i< g$, and each in $(V(\mathbb{C})\setminus \{v_{i-1}$, $v_{i}$, $v_{i+1}\})\cup V(\mathbb{P}-v_{g+l})$ is a $p$-dominator (if $i=1$, then $v_{i-1}=v_{g}$).
\end{theorem}

\begin{proof}
Denote by $A$ the set of vertices of $\mathbb{C}$ and the pendant vertices attached to $\mathbb{C}$. Let $\parallel A \parallel=z$, and let $A^{'}=V(\mathcal {G})\setminus A$. Then $\gamma(\mathcal {G})\leq \gamma(\mathcal {G}[A])+\gamma(\mathcal {G}[A^{'}])$. Note that $A^{'}=\emptyset$, or $\mathcal {G}[A^{'}]$ is connected with at least $2$ vertices. Suppose $f\geq 4$.

(i) $f=g$. Then $z-f=0$. This means that there is no $p$-dominator on $\mathbb{C}$. So, $\mathcal {G}[A^{'}]$ is connected with at least $2$ vertices. Thus, if $f\geq 9$, by Lemma \ref{le0,09}, then $\gamma(\mathcal {G})\leq \lceil\frac{f}{3}\rceil+\gamma(\mathcal {G}[A^{'}])\leq \frac{n-f}{2}+\frac{f+2}{3}<\frac{n-1}{2}$. Therefore $f\leq 7$.

Note that $g$ is odd and $g=f$ now. Thus if $\gamma(\mathcal {G}[A^{'}])<\frac{n-f}{2}$, then $\gamma(\mathcal {G})\leq \lceil\frac{f}{3}\rceil+\gamma(\mathcal {G}[A^{'}])<\frac{n-1}{2}$. Hence, it follows that $\gamma(\mathcal {G}[A^{'}])=\frac{n-f}{2}$. Combined with Lemma \ref{le0,10}, it follows that $\mathcal {G}[A^{'}]=P_{\frac{n-f}{2}}\circ K_{1}$. Here, suppose $P_{\frac{n-f}{2}}=v_{a_{1}}v_{a_{2}}\cdots v_{a_{t}}$ with $t=\frac{n-f}{2}$, and suppose $v_{\tau_{1}}$ is the unique pendant vertex attached to $v_{a_{1}}$. By Lemma \ref{le0,08}, $V(P_{\frac{n-f}{2}})$ is a minimal dominating set of $\mathcal {G}[A^{'}]$.

Assume that $f=7$. Note that $\mathcal {G}$ is a $\mathcal {F}_{g, l}$-graph. If $\mathcal {G}=\mathbb{C}+v_{g}v_{a_{1}}+\mathcal {G}[A^{'}]$, then $V(P_{\frac{n-f}{2}})\cup \{v_{2}, v_{5}\}$ is a dominating set of $\mathcal {G}$; if $\mathcal {G}=\mathbb{C}+v_{g}v_{\tau_{1}}+\mathcal {G}[A^{'}]$, then $(V(P_{\frac{n-f}{2}})\backslash \{v_{a_{1}}\})\cup \{v_{2}, v_{5}, v_{\tau_{1}}\}$ is a dominating set of $\mathcal {G}$. This implies that $\gamma(\mathcal {G})\leq \frac{n-7}{2}+2<\frac{n-1}{2}$ which contradicts $\gamma(\mathcal {G}) =\frac{n-1}{2}$. Thus, it follows that $g=5$.

(ii) $f\neq g$. Note that there is no the case that $z-f=1$. Then $z-f\geq2$.
By Lemma \ref{le03,07}, $\gamma(\mathcal {G}[A])\leq \gamma(\mathcal {G}^{\ast}[A])=g-f+\lceil\frac{f-2}{3}\rceil\leq \frac{z-f}{2}+\lceil\frac{f-2}{3}\rceil$, where $\mathcal {G}^{\ast}[A]$ is a sunlike graph with vertex set $A$, $\mathbb{C}$ contained in it and $g-f$ $p$-dominators $v_{1}$, $v_{2}$, $\ldots$, $v_{g-f}$ (defined as $\mathcal {G}^{\ast}$ in Lemma \ref{le03,07}). Thus, if $f\geq 4$, then $\gamma(\mathcal {G})\leq \frac{z-f}{2}+\lceil\frac{f-2}{3}\rceil+\gamma(\mathcal {G}[A^{'}])\leq \frac{n-f}{2}+\lceil\frac{f-2}{3}\rceil\leq \frac{n-f}{2}+\frac{f}{3}<\frac{n-1}{2}$. This contradicts that $\gamma(\mathcal {G}) =\frac{n-1}{2}$. Consequently, $f\leq 3$.

Suppose $f=2$ and suppose that $v_{j}$, $v_{k}$ of $\mathbb{C}$ are the exact $2$ vertices such that neither of them is $p$-dominator. Note that by Lemma \ref{le0,08}, there is a minimal dominating set $D$ of $\mathcal {G}-v_{j}-v_{k}$ which contains all $p$-dominators but no any pendant vertex. Note that the vertices of $\mathbb{C}$ other than $v_{j}$, $v_{k}$ are all $p$-dominators in both $\mathcal {G}-v_{j}-v_{k}$ and $\mathcal {G}$. Thus, each of $v_{j}$, $v_{k}$ is adjacent to at least one $p$-dominator on $\mathbb{C}$. So, $D$ is also a dominating set of $\mathcal {G}$. Note that there is no isolated vertex in $\mathcal {G}-v_{j}-v_{k}$. Then $\gamma(\mathcal {G}-v_{j}-v_{k})\leq\frac{n-2}{2}$, and then $\gamma(\mathcal {G})\leq\frac{n-2}{2}$, which contradicts $\gamma(\mathcal {G}) =\frac{n-1}{2}$. Then (ii) follows.

(iii) Suppose $v_{a}$, $v_{b}$, $v_{c}$ are the exact $3$ vertices of $\mathbb{C}$ such that none of them is $p$-dominator. If the $3$ vertices $v_{a}$, $v_{b}$, $v_{c}$ are not
consecutive, then each of them can be dominated by its adjacent $p$-dominator. Note that by Lemma \ref{le0,08}, there are a minimal dominating set $D$ of $\mathcal {G}-v_{a}-v_{b}-v_{c}$ which contains all $p$-dominators but no any pendant vertex. Thus such $D$ is also a dominating set of $\mathcal {G}$. Note that there is no isolated vertex in $\mathcal {G}-v_{a}-v_{b}-v_{c}$. So, $\gamma(\mathcal {G})\leq \parallel D \parallel=\gamma (\mathcal {G}-v_{a}-v_{b}-v_{c})\leq \frac{n-3}{2}$, which contradicts $\gamma(\mathcal {G}) =\frac{n-1}{2}$. Therefore, the $3$ vertices $v_{a}$, $v_{b}$, $v_{c}$ are
consecutive.

Suppose that the $3$ vertices are $v_{i-1}$, $v_{i}$, $v_{i+1}$ for some $1\leq i\leq g$ (here, if $i=g$, we let $v_{i+1}=v_{1}$; if $i=1$, we let $v_{i-1}=v_{g}$). Let $H=\mathcal {G}-v_{i-1}-v_{i}-v_{i+1}$. Note that there is no isolated vertex in $H$. Thus, $\gamma(H)\leq\frac{n-3}{2}$. Next, we claim that $\gamma(H)=\frac{n-3}{2}$.

{\bf Claim 1} $\gamma(H)=\frac{n-3}{2}$. Otherwise, suppose $\gamma(H)<\frac{n-3}{2}$, and suppose $D$ is a minimal dominating set of $H$. Then $D\cup\{v_{i}\}$ is a dominating set $D$ of $\mathcal {G}$. Thus, $\mathcal {G}< 1+\frac{n-3}{2}<\frac{n-1}{2}$, which contradicts $\gamma(\mathcal {G}) =\frac{n-1}{2}$. Then the claim holds.

By Lemma \ref{le0,10}, $H=\mathcal {L}\circ K_{1}$ for some acyclic graph $\mathcal {L}$ of order $\frac{n-3}{2}$.

{\bf Claim 2} For any minimal dominating set $D$ of $H$, in $\mathcal {G}$, at least one of $v_{i-1}$, $v_{i}$, $v_{i+1}$ can not be dominated by $D$. Otherwise, $D$ is a dominating set of $\mathcal {G}$ too. Hence, $\gamma(\mathcal {G})\leq\frac{n-3}{2}$, which contradicts $\gamma(\mathcal {G}) =\frac{n-1}{2}$. Then the claim holds.

If $i=g$, then let $H=H_{1}\cup H_{2}$, where $H_{1}=\mathcal {G}[A]-v_{g-1}-v_{g}-v_{1}$, $H_{2}=\mathcal {G}[A^{'}]=P_{\frac{n-z}{2}}\circ K_{1}$ (if $n=z$, then $H_{2}$ is empty). Here, suppose $P_{\frac{n-z}{2}}=v_{a_{1}}v_{a_{2}}\cdots v_{a_{t}}$ with $t=\frac{n-z}{2}$, and suppose $v_{\tau_{1}}$ is the unique pendant vertex attached to $v_{a_{1}}$. Thus there are two possible cases for $G$, i.e., $\mathcal {G}=\mathcal {G}[A]+v_{g}v_{a_{1}}+H_{2}$ or $\mathcal {G}=\mathcal {G}[A]+v_{g}v_{\tau_{1}}+H_{2}$. Let $\mathcal {Z}=(\mathbb{C}\setminus \{v_{g-1}, v_{g}, v_{1}\})\cup V(P_{\frac{n-z}{2}})$. Note that the vertices in $\mathcal {Z}$ are all $p$-dominators in $\mathcal {G}$. If $\mathcal {G}=\mathcal {G}[A]+v_{g}v_{a_{1}}+H_{2}$, then $\mathcal {Z}$ is also a dominating set of $\mathcal {G}$; if $\mathcal {G}=\mathcal {G}[A]+v_{g}v_{\tau_{1}}+H_{2}$, then $(\mathcal {Z}\setminus \{v_{a_{1}}\})\cup\{v_{\tau_{1}}\}$ is a dominating set of $\mathcal {G}$. Thus it follows that $\gamma(\mathcal {G})\leq \frac{n-3}{2}<\frac{n-1}{2}$ which contradicts $\gamma(\mathcal {G}) =\frac{n-1}{2}$. This implies $i\neq g$.

If $i\neq 1, g-1$, then $H$ is connected. Let $\mathcal {Z}=(V(\mathbb{C})\setminus \{v_{i-1}$, $v_{i}$, $v_{i+1}\})\cup V(\mathbb{P}-v_{g+l})$, where $\mathbb{P}=v_{g}v_{g+1}\cdots v_{g+l}$. Then each vertex in $\mathcal {Z}$ is a $p$-dominator in $\mathcal {G}$.

If $i= 1$, then let $H=H_{1}\cup H_{2}$, where $H_{1}=\mathcal {G}[A]-v_{g}-v_{1}-v_{2}$, $H_{2}=\mathcal {G}[A^{'}]=P_{\frac{n-z}{2}}\circ K_{1}$ (if $n=z$, then $H_{2}$ is empty). Here, suppose $P_{\frac{n-z}{2}}=v_{a_{1}}v_{a_{2}}\cdots v_{a_{t}}$ with $t=\frac{n-z}{2}$, and suppose $v_{\tau_{1}}$ is the unique pendant vertex attached to $v_{a_{1}}$. Thus there are two possible cases for $G$, i.e., $\mathcal {G}=\mathcal {G}[A]+v_{g}v_{a_{1}}+H_{2}$ or $\mathcal {G}=\mathcal {G}[A]+v_{g}v_{\tau_{1}}+H_{2}$. We say that $\mathcal {G}\neq\mathcal {G}[A]+v_{g}v_{\tau_{1}}+H_{2}$. Otherwise, suppose $\mathcal {G}=\mathcal {G}[A]+v_{g}v_{\tau_{1}}+H_{2}$. Note that $n-z$ is even now and $\mathcal {G}-\{v_{2}, v_{1}, v_{g}, v_{a_{1}}, v_{\tau_{1}}\}$ has no isolated vertex. Then for $\mathcal {G}-\{v_{2}, v_{1}, v_{g}, v_{a_{1}}, v_{\tau_{1}}\}$, it has a dominating set $\mathbb{D}$ with $\parallel\mathbb{D}\parallel\leq\frac{n-5}{2}$. Then $\mathbb{D}\cup \{v_{1}, v_{\tau_{1}}\}$ is a dominating set of $\mathcal {G}$, which contradicts $\gamma(\mathcal {G}) =\frac{n-1}{2}$. This implies that $\mathcal {G}=\mathcal {G}[A]+v_{g}v_{a_{1}}+H_{2}$. It follows that each one in $(V(\mathbb{C})\setminus \{v_{g}$, $v_{1}$, $v_{2}\})\cup V(\mathbb{P}-v_{g+l})$ is a $p$-dominator. Similarly, for $i= g-1$, we get that each one in $(V(\mathbb{C})\setminus \{v_{g-2}$, $v_{g-1}$, $v_{g}\})\cup V(\mathbb{P}-v_{g+l})$ is a $p$-dominator.
Then (iii) follows.
 \ \ \ \ \ $\Box$
\end{proof}

\section{The $q_{min}$ among uncyclic graphs}

\begin{lemma}{\bf \cite{YMCG}}\label{le04,01} 
Let $G$ be a nonbipartite unicyclic graph of order $n$ and with the odd cycle $\mathcal {C}=v_1v_2\cdots v_gv_1$ in it. There is a
unit eigenvector $X=(\,x_1$, $x_2$, $\ldots$, $x_g$, $x_{g+1}$, $x_{g+2}$, $\ldots$, $x_{n-1}$, $x_{n}\,)^T$ corresponding to $q_{min}(G)$, in which suppose $|x_1|=\min\{|x_1|$, $|x_2|$, $\ldots$, $|x_g|\}$,
$|x_s|=\max\{|x_1|$, $|x_2|$, $\ldots$, $|x_g|\}$ where $s\geq 2$, satisfying that

$\mathrm{(i)}$ $|x_1|< |x_s|$;

$\mathrm{(ii)}$ $|x_1|=0$ if and only if $x_{g}=-x_{2}\neq 0$; if $|x_1|=0$ and $x_ix_{i+1}\neq 0$ for some $1\leq i\leq g-1$, then $x_ix_{i+1}< 0$; moreover, if $x_{j}\neq 0$, then $sgn(x_{j})=(-1)^{d_{istH}(v_{1}, v_{j})}$ where $H=G-v_1v_g$.

$\mathrm{(iii)}$ if $|x_1|>0$, then

$\mathrm{(1)}$  if $3\leq s\leq g-1$, then $|x_2|<\cdots<|x_{s-2}|< |x_{s-1}|\leq |x_s|$ and $|x_g|<|x_{g-1}|<\cdots<|x_{s+2}|< |x_{s+1}|\leq |x_s|$;

$\mathrm{(2)}$ if $|x_2|> |x_g|$, then $x_1x_{g}> 0$; for $1\leq i\leq g-1$, $x_ix_{i+1}< 0$; $|x_1|\leq |x_g|$;

$\mathrm{(3)}$ if $|x_2|<|x_g|$,
then $x_1x_{2}> 0$; for $2\leq i\leq g-1$, $x_ix_{i+1}< 0$; $x_gx_{1}< 0$; $|x_1|\leq|x_2|$;

$\mathrm{(4)}$ if $|x_2|=|x_g|$, then $|x_1|\leq|x_2|$, and exactly one of $x_1x_{g}> 0$ and $x_1x_{2}> 0$ holds, where

\ \ \ \ \ $\mathrm{(4.1)}$ if $x_1x_{g}> 0$, then for $1\leq i\leq g-1$, $x_ix_{i+1}< 0$;

\ \ \ \ \ $\mathrm{(4.2)}$ if $x_1x_{2}> 0$, then $x_ix_{i+1}< 0$ for $2\leq i\leq g-1$ and $x_gx_{1}< 0$;

$\mathrm{(5)}$ at least one of  $|x_{s+1}|$ and $|x_{s-1}|$ is less than $|x_{s}|$.
\end{lemma}

\begin{lemma}{\bf \cite{YMCG}}\label{th04,02} 
If $\mathcal {G}$ is a nonbipartite $\mathcal {F}^{\circ}_{g, l}$-graph with $g\geq 5$, $n\geq g+1$, then there is a graph $\mathbb{H}$ with girth $3$ and order $n$ such that $\gamma(\mathcal {G})\leq\gamma(\mathbb{H})$ and $q_{min}(\mathbb{H})<q_{min}(\mathcal {G})$.
\end{lemma}

\begin{lemma}{\bf \cite{YMCG}}\label{le04,03} 
Suppose that $G$ is a nonbipartite $\mathcal {F}_{3, l}$-graph of order $n$ where $\mathbb{C}=v_{1}v_{2}v_{3}v_{1}$. $X=(\,x_{1}$, $x_{2}$, $\ldots$, $x_{n}\,)^T$ is a unit eigenvector corresponding to $q_{min}(G)$. Then $|x_{3}|=\max\{|x_{1}|$, $|x_{2}|$, $|x_{3}|\}$.
\end{lemma}

\begin{theorem}\label{th04,04} 
Among all nonbipartite unicyclic graphs of order $n\geq 5$ with girth $3$ and domination number at least $\frac{n+1}{3}<\gamma\leq \frac{n}{2}$,
if $\gamma=\frac{n-1}{2}$, the $q_{min}$ attains the minimum uniquely at $\mathscr{H}_{3,\frac{n-3}{2}}$.
\end{theorem}

\begin{proof}
The result follows from Lemmas \ref{le0,04}, \ref{th03,02}, \ref{th03,05}, \ref{le04,03} and Theorem \ref{le03,08}
 \ \ \ \ \ $\Box$
\end{proof}

Let $\mathcal {K}=\{G|\,\ G\ \mathrm{be}\ \mathrm{a}\ \mathrm{nonbipartite}\ \mathcal {F}^{\circ}_{g, l}$-graph of order $n\geq 4$ and domination number at least $\frac{n+1}{3}<\gamma\leq \frac{n}{2}$, where $g$ is any odd number at least $3$ and $l$ is any positive integral number$\}$ and $q_{\mathcal {K}}=\min\{q_{min}(G)|\,\ G\in \mathcal {K}\}$.

\begin{lemma}{\bf \cite{YMCG}}\label{th04,05} 
\

$\mathrm{(i)}$ If $n=4$, the $q_{\mathcal {K}}$ attains uniquely at $\mathscr{H}_{3,1}$;

$\mathrm{(ii)}$ If $n\geq 5$ and $n-2\gamma\geq 2$, then the least $q_{\mathcal {K}}> q_{min}(\mathscr{H}_{3,\alpha})$ where $\alpha\leq\frac{n-3}{2}$ is the least integer such that $\lceil\frac{n-2\alpha-2}{3}\rceil+\alpha=\gamma$.
\end{lemma}

\begin{lemma}\label{le04,09} 
For  a nonbipartite $\mathcal {F}_{g, l}$-graph graph $G$ of order $n\geq 5$ and with $g =5$, there exists a graph $\mathbb{H}$ such that $g(\mathbb{H})=3$, $\gamma(G)\leq \gamma(\mathbb{H})$ and $q_{min}(\mathbb{H})< q_{min}(G)$.
\end{lemma}

\begin{proof}
If $n=5$, then $G=C_{5}$. And then the result follows from Lemma \ref{th2,18}. Next we consider the case that $n\geq 6$. By Lemma \ref{le0,06}, we get that $q_{min}(G)<1$.

{\bf Case 1} There is no $p$-dominator on $\mathbb{C}$. Then $G$ is like $G_{1}$ (see $G_{1}$ in Fig. 4.1). By Lemma \ref{le0,05}, there is a unit eigenvector $X=(\,x_1$, $x_2$, $\ldots$, $x_k$, $x_{k+1}$, $x_{k+2}$, $\ldots$, $x_{n-1}$, $x_{n}\,)^T$ corresponding to $q_{min}(G)$ such that $|x_5|=\max\{|x_1|$, $|x_2|$, $|x_3|$, $|x_4|$, $|x_5|\}> 0$, and $x_{1}=x_{4}$, $x_{2}=x_{3}$. By Lemma \ref{le04,01}, we get that $|x_2|>0$, $|x_2|< |x_1|$ and $x_2 x_1< 0$. Let $\mathbb{H}=G-v_{3}v_{4}+v_{3}v_{1}$. By Lemma \ref{le0,04}, we get that $q_{min}(\mathbb{H})< q_{min}(G)$. Let $B_{1}=\mathbb{H}[v_{1}$, $v_{2}$, $v_{3}]$, $B_{2}=\mathbb{H}- \{v_{1}$, $v_{2}$, $v_{3}\}$. As Lemma \ref{le03,07}, we can get a minimal dominating set $D$ of $\mathbb{H}$,
which contains all $p$-dominators but no any pendant vertex and no $v_{1}$, such that $D=\{v_{2}\}\cup D_{2}$, where $\{v_{2}\}$ is a dominating set of $B_{1}$, $D_{2}$ is a dominating set of $B_{2}$. Note that $D$ is also a dominating set of $G$. So, $\gamma(G)\leq \gamma(\mathbb{H})$.

\setlength{\unitlength}{0.6pt}
\begin{center}
\begin{picture}(636,599)
\put(38,577){\circle*{4}}
\put(19,541){\circle*{4}}
\qbezier(38,577)(28,559)(19,541)
\put(79,546){\circle*{4}}
\qbezier(38,577)(58,562)(79,546)
\put(19,507){\circle*{4}}
\qbezier(19,541)(19,524)(19,507)
\put(79,507){\circle*{4}}
\qbezier(79,546)(79,527)(79,507)
\qbezier(19,507)(49,507)(79,507)
\put(86,577){\circle*{4}}
\qbezier(38,577)(62,577)(86,577)
\put(127,577){\circle*{4}}
\qbezier(86,577)(106,577)(127,577)
\put(86,598){\circle*{4}}
\qbezier(86,577)(86,588)(86,598)
\put(127,598){\circle*{4}}
\qbezier(127,577)(127,588)(127,598)
\put(1,540){$v_{1}$}
\put(1,500){$v_{2}$}
\put(84,502){$v_{3}$}
\put(84,543){$v_{4}$}
\put(20,582){$v_{5}$}
\put(78,565){$v_{5+s}$}
\put(36,473){$G_{1}$}
\put(253,577){\circle*{4}}
\put(270,544){\circle*{4}}
\qbezier(253,577)(261,561)(270,544)
\put(209,541){\circle*{4}}
\qbezier(253,577)(231,559)(209,541)
\put(209,507){\circle*{4}}
\qbezier(209,541)(209,524)(209,507)
\put(270,507){\circle*{4}}
\qbezier(270,544)(270,526)(270,507)
\qbezier(209,507)(239,507)(270,507)
\put(163,577){\circle*{4}}
\qbezier(253,577)(208,577)(163,577)
\put(163,598){\circle*{4}}
\qbezier(163,577)(163,588)(163,598)
\put(205,577){\circle*{4}}
\put(205,598){\circle*{4}}
\qbezier(205,577)(205,588)(205,598)
\put(253,598){\circle*{4}}
\qbezier(253,577)(253,588)(253,598)
\put(191,541){$v_{1}$}
\put(190,504){$v_{2}$}
\put(274,504){$v_{3}$}
\put(274,544){$v_{4}$}
\put(258,576){$v_{5}$}
\put(230,473){$G_{2}$}
\put(363,577){\circle*{4}}
\put(349,549){\circle*{4}}
\qbezier(363,577)(356,563)(349,549)
\put(405,549){\circle*{4}}
\qbezier(363,577)(384,563)(405,549)
\put(349,508){\circle*{4}}
\qbezier(349,549)(349,529)(349,508)
\put(405,508){\circle*{4}}
\qbezier(405,549)(405,529)(405,508)
\qbezier(349,508)(377,508)(405,508)
\put(456,577){\circle*{4}}
\qbezier(363,577)(409,577)(456,577)
\put(456,598){\circle*{4}}
\qbezier(456,577)(456,588)(456,598)
\put(414,577){\circle*{4}}
\put(414,598){\circle*{4}}
\qbezier(414,577)(414,588)(414,598)
\put(442,549){\circle*{4}}
\qbezier(405,549)(423,549)(442,549)
\put(331,547){$v_{1}$}
\put(330,504){$v_{2}$}
\put(409,504){$v_{3}$}
\put(388,541){$v_{4}$}
\put(353,585){$v_{5}$}
\put(402,566){$v_{5+s}$}
\put(370,473){$G_{3}$}
\put(592,577){\circle*{4}}
\put(608,548){\circle*{4}}
\qbezier(592,577)(600,563)(608,548)
\put(548,546){\circle*{4}}
\qbezier(592,577)(570,562)(548,546)
\put(548,505){\circle*{4}}
\qbezier(548,546)(548,526)(548,505)
\put(608,505){\circle*{4}}
\qbezier(608,548)(608,527)(608,505)
\qbezier(548,505)(578,505)(608,505)
\put(497,577){\circle*{4}}
\qbezier(497,577)(544,577)(592,577)
\put(497,597){\circle*{4}}
\qbezier(497,577)(497,587)(497,597)
\put(541,577){\circle*{4}}
\put(541,598){\circle*{4}}
\qbezier(541,577)(541,588)(541,598)
\put(512,505){\circle*{4}}
\qbezier(548,505)(530,505)(512,505)
\put(530,545){$v_{1}$}
\put(550,510){$v_{2}$}
\put(613,502){$v_{3}$}
\put(614,548){$v_{4}$}
\put(592,582){$v_{5}$}
\put(523,565){$v_{s+5}$}
\put(568,473){$G_{4}$}
\put(35,415){\circle*{4}}
\put(18,381){\circle*{4}}
\qbezier(35,415)(26,398)(18,381)
\put(73,385){\circle*{4}}
\qbezier(35,415)(54,400)(73,385)
\put(18,343){\circle*{4}}
\qbezier(18,381)(18,362)(18,343)
\put(73,343){\circle*{4}}
\qbezier(73,385)(73,364)(73,343)
\qbezier(18,343)(45,343)(73,343)
\put(101,415){\circle*{4}}
\qbezier(35,415)(68,415)(101,415)
\put(101,441){\circle*{4}}
\qbezier(101,415)(101,428)(101,441)
\put(71,415){\circle*{4}}
\put(71,440){\circle*{4}}
\qbezier(71,415)(71,428)(71,440)
\put(35,440){\circle*{4}}
\qbezier(35,415)(35,428)(35,440)
\put(0,381){$v_{1}$}
\put(8,331){$v_{2}$}
\put(77,338){$v_{3}$}
\put(55,380){$v_{4}$}
\put(16,417){$v_{5}$}
\put(102,385){\circle*{4}}
\qbezier(73,385)(87,385)(102,385)
\put(40,306){$G_{5}$}
\put(190,415){\circle*{4}}
\put(207,383){\circle*{4}}
\qbezier(190,415)(198,399)(207,383)
\put(151,381){\circle*{4}}
\qbezier(190,415)(170,398)(151,381)
\put(151,343){\circle*{4}}
\qbezier(151,381)(151,362)(151,343)
\put(207,343){\circle*{4}}
\qbezier(151,343)(179,343)(207,343)
\qbezier(207,383)(207,363)(207,343)
\put(125,415){\circle*{4}}
\qbezier(190,415)(157,415)(125,415)
\put(125,440){\circle*{4}}
\qbezier(125,415)(125,428)(125,440)
\put(153,415){\circle*{4}}
\put(153,441){\circle*{4}}
\qbezier(153,415)(153,428)(153,441)
\put(190,441){\circle*{4}}
\qbezier(190,415)(190,428)(190,441)
\put(123,343){\circle*{4}}
\qbezier(151,343)(137,343)(123,343)
\put(155,375){$v_{1}$}
\put(140,332){$v_{2}$}
\put(211,340){$v_{3}$}
\put(210,382){$v_{4}$}
\put(194,415){$v_{5}$}
\put(160,306){$G_{6}$}
\put(276,413){\circle*{4}}
\put(261,386){\circle*{4}}
\qbezier(276,413)(268,400)(261,386)
\put(261,341){\circle*{4}}
\qbezier(261,386)(261,364)(261,341)
\put(314,383){\circle*{4}}
\qbezier(276,413)(295,398)(314,383)
\put(314,341){\circle*{4}}
\qbezier(314,383)(314,362)(314,341)
\qbezier(261,341)(287,341)(314,341)
\put(344,413){\circle*{4}}
\qbezier(276,413)(310,413)(344,413)
\put(344,440){\circle*{4}}
\qbezier(344,413)(344,427)(344,440)
\put(314,413){\circle*{4}}
\put(314,440){\circle*{4}}
\qbezier(314,413)(314,427)(314,440)
\put(342,383){\circle*{4}}
\qbezier(314,383)(328,383)(342,383)
\put(343,341){\circle*{4}}
\qbezier(314,341)(328,341)(343,341)
\put(242,383){$v_{1}$}
\put(242,338){$v_{2}$}
\put(312,329){$v_{3}$}
\put(296,378){$v_{4}$}
\put(258,416){$v_{5}$}
\put(301,402){$v_{s+5}$}
\put(283,306){$G_{7}$}
\put(441,414){\circle*{4}}
\put(402,381){\circle*{4}}
\qbezier(441,414)(421,398)(402,381)
\put(456,383){\circle*{4}}
\qbezier(441,414)(448,399)(456,383)
\put(402,343){\circle*{4}}
\qbezier(402,381)(402,362)(402,343)
\put(456,343){\circle*{4}}
\qbezier(456,383)(456,363)(456,343)
\qbezier(402,343)(429,343)(456,343)
\put(372,414){\circle*{4}}
\qbezier(441,414)(406,414)(372,414)
\put(372,441){\circle*{4}}
\qbezier(372,414)(372,428)(372,441)
\put(403,414){\circle*{4}}
\put(403,442){\circle*{4}}
\qbezier(403,414)(403,428)(403,442)
\put(375,381){\circle*{4}}
\qbezier(402,381)(388,381)(375,381)
\put(484,343){\circle*{4}}
\qbezier(456,343)(470,343)(484,343)
\put(407,375){$v_{1}$}
\put(384,336){$v_{2}$}
\put(449,332){$v_{3}$}
\put(460,380){$v_{4}$}
\put(440,418){$v_{5}$}
\put(390,403){$v_{s+5}$}
\put(423,306){$G_{8}$}
\put(554,415){\circle*{4}}
\put(538,383){\circle*{4}}
\qbezier(554,415)(546,399)(538,383)
\put(590,384){\circle*{4}}
\qbezier(554,415)(572,400)(590,384)
\put(538,343){\circle*{4}}
\qbezier(538,383)(538,363)(538,343)
\put(590,343){\circle*{4}}
\qbezier(590,384)(590,364)(590,343)
\qbezier(538,343)(564,343)(590,343)
\put(627,415){\circle*{4}}
\qbezier(554,415)(590,415)(627,415)
\put(627,444){\circle*{4}}
\qbezier(627,415)(627,430)(627,444)
\put(594,415){\circle*{4}}
\put(594,443){\circle*{4}}
\qbezier(594,415)(594,429)(594,443)
\put(512,383){\circle*{4}}
\qbezier(538,383)(525,383)(512,383)
\put(617,384){\circle*{4}}
\qbezier(590,384)(603,384)(617,384)
\put(520,371){$v_{1}$}
\put(519,341){$v_{2}$}
\put(594,343){$v_{3}$}
\put(572,378){$v_{4}$}
\put(536,418){$v_{5}$}
\put(585,404){$v_{s+5}$}
\put(556,306){$G_{9}$}
\put(72,253){\circle*{4}}
\put(38,219){\circle*{4}}
\qbezier(72,253)(55,236)(38,219)
\put(92,222){\circle*{4}}
\qbezier(72,253)(82,238)(92,222)
\put(38,181){\circle*{4}}
\qbezier(38,219)(38,200)(38,181)
\put(92,181){\circle*{4}}
\qbezier(92,222)(92,202)(92,181)
\qbezier(38,181)(65,181)(92,181)
\put(8,253){\circle*{4}}
\qbezier(72,253)(40,253)(8,253)
\put(8,282){\circle*{4}}
\qbezier(8,253)(8,268)(8,282)
\put(37,253){\circle*{4}}
\put(37,282){\circle*{4}}
\qbezier(37,253)(37,268)(37,282)
\put(14,181){\circle*{4}}
\qbezier(38,181)(26,181)(14,181)
\put(119,181){\circle*{4}}
\qbezier(92,181)(105,181)(119,181)
\put(20,219){$v_{1}$}
\put(30,169){$v_{2}$}
\put(83,170){$v_{3}$}
\put(96,219){$v_{4}$}
\put(75,255){$v_{5}$}
\put(46,144){$G_{10}$}
\put(22,242){$v_{s+5}$}
\put(167,253){\circle*{4}}
\put(148,221){\circle*{4}}
\qbezier(167,253)(157,237)(148,221)
\put(148,181){\circle*{4}}
\qbezier(148,221)(148,201)(148,181)
\put(200,220){\circle*{4}}
\qbezier(167,253)(183,237)(200,220)
\put(200,181){\circle*{4}}
\qbezier(200,220)(200,201)(200,181)
\qbezier(148,181)(174,181)(200,181)
\put(236,253){\circle*{4}}
\qbezier(167,253)(201,253)(236,253)
\put(236,282){\circle*{4}}
\qbezier(236,253)(236,268)(236,282)
\put(208,253){\circle*{4}}
\put(208,283){\circle*{4}}
\qbezier(208,253)(208,268)(208,283)
\put(135,253){\circle*{4}}
\qbezier(167,253)(151,253)(135,253)
\put(229,220){\circle*{4}}
\qbezier(200,220)(214,220)(229,220)
\put(230,181){\circle*{4}}
\qbezier(200,181)(215,181)(230,181)
\put(131,219){$v_{1}$}
\put(141,168){$v_{2}$}
\put(198,168){$v_{3}$}
\put(182,216){$v_{4}$}
\put(159,260){$v_{5}$}
\put(163,144){$G_{11}$}
\put(335,253){\circle*{4}}
\put(295,219){\circle*{4}}
\qbezier(335,253)(315,236)(295,219)
\put(295,181){\circle*{4}}
\qbezier(295,219)(295,200)(295,181)
\put(349,220){\circle*{4}}
\qbezier(335,253)(342,237)(349,220)
\put(349,181){\circle*{4}}
\qbezier(349,220)(349,201)(349,181)
\qbezier(295,181)(322,181)(349,181)
\put(261,253){\circle*{4}}
\qbezier(335,253)(298,253)(261,253)
\put(261,280){\circle*{4}}
\qbezier(261,253)(261,267)(261,280)
\put(296,253){\circle*{4}}
\put(296,281){\circle*{4}}
\qbezier(296,253)(296,267)(296,281)
\put(368,253){\circle*{4}}
\qbezier(335,253)(351,253)(368,253)
\put(269,219){\circle*{4}}
\qbezier(295,219)(282,219)(269,219)
\put(376,220){\circle*{4}}
\qbezier(349,220)(362,220)(376,220)
\put(282,225){$v_{1}$}
\put(276,178){$v_{2}$}
\put(353,177){$v_{3}$}
\put(331,219){$v_{4}$}
\put(329,258){$v_{5}$}
\put(310,144){$G_{12}$}
\put(428,255){\circle*{4}}
\put(411,221){\circle*{4}}
\qbezier(428,255)(419,238)(411,221)
\put(467,222){\circle*{4}}
\qbezier(428,255)(447,239)(467,222)
\put(411,181){\circle*{4}}
\qbezier(411,221)(411,201)(411,181)
\put(467,181){\circle*{4}}
\qbezier(467,222)(467,202)(467,181)
\qbezier(411,181)(439,181)(467,181)
\put(498,255){\circle*{4}}
\qbezier(428,255)(463,255)(498,255)
\put(469,255){\circle*{4}}
\put(469,278){\circle*{4}}
\qbezier(469,255)(469,267)(469,278)
\put(498,279){\circle*{4}}
\qbezier(498,255)(498,267)(498,279)
\put(496,222){\circle*{4}}
\qbezier(467,222)(481,222)(496,222)
\put(496,181){\circle*{4}}
\qbezier(467,181)(481,181)(496,181)
\put(399,245){\circle*{4}}
\qbezier(411,221)(405,233)(399,245)
\put(424,144){$G_{13}$}
\put(414,212){$v_{1}$}
\put(392,175){$v_{2}$}
\put(463,169){$v_{3}$}
\put(450,216){$v_{4}$}
\put(421,261){$v_{5}$}
\put(458,244){$v_{5+s}$}
\put(590,255){\circle*{4}}
\put(551,221){\circle*{4}}
\qbezier(590,255)(570,238)(551,221)
\put(606,222){\circle*{4}}
\qbezier(590,255)(598,239)(606,222)
\put(551,181){\circle*{4}}
\qbezier(551,221)(551,201)(551,181)
\put(606,181){\circle*{4}}
\qbezier(606,222)(606,202)(606,181)
\qbezier(551,181)(578,181)(606,181)
\put(524,181){\circle*{4}}
\qbezier(551,181)(537,181)(524,181)
\put(524,221){\circle*{4}}
\qbezier(551,221)(537,221)(524,221)
\put(636,181){\circle*{4}}
\qbezier(606,181)(621,181)(636,181)
\put(522,255){\circle*{4}}
\qbezier(590,255)(556,255)(522,255)
\put(522,279){\circle*{4}}
\qbezier(522,255)(522,267)(522,279)
\put(555,255){\circle*{4}}
\put(555,280){\circle*{4}}
\qbezier(555,255)(555,268)(555,280)
\put(554,214){$v_{1}$}
\put(543,169){$v_{2}$}
\put(603,169){$v_{3}$}
\put(610,219){$v_{4}$}
\put(587,260){$v_{5}$}
\put(540,244){$v_{5+s}$}
\put(567,144){$G_{14}$}
\put(47,97){\circle*{4}}
\put(25,64){\circle*{4}}
\qbezier(47,97)(36,81)(25,64)
\put(80,66){\circle*{4}}
\qbezier(47,97)(63,82)(80,66)
\put(25,23){\circle*{4}}
\qbezier(25,64)(25,44)(25,23)
\put(80,23){\circle*{4}}
\qbezier(80,66)(80,45)(80,23)
\qbezier(25,23)(52,23)(80,23)
\put(118,97){\circle*{4}}
\qbezier(47,97)(82,97)(118,97)
\put(118,122){\circle*{4}}
\qbezier(118,97)(118,110)(118,122)
\put(90,97){\circle*{4}}
\put(90,122){\circle*{4}}
\qbezier(90,97)(90,110)(90,122)
\put(17,97){\circle*{4}}
\qbezier(47,97)(32,97)(17,97)
\put(111,66){\circle*{4}}
\qbezier(80,66)(95,66)(111,66)
\put(111,23){\circle*{4}}
\qbezier(80,23)(95,23)(111,23)
\put(9,34){\circle*{4}}
\qbezier(25,23)(17,29)(9,34)
\put(29,58){$v_{1}$}
\put(18,12){$v_{2}$}
\put(78,12){$v_{3}$}
\put(62,60){$v_{4}$}
\put(42,103){$v_{5}$}
\put(43,-17){$G_{15}$}
\put(176,97){\circle*{4}}
\put(154,65){\circle*{4}}
\qbezier(176,97)(165,81)(154,65)
\put(154,23){\circle*{4}}
\qbezier(154,65)(154,44)(154,23)
\put(205,66){\circle*{4}}
\qbezier(176,97)(190,82)(205,66)
\put(205,23){\circle*{4}}
\qbezier(205,66)(205,45)(205,23)
\qbezier(154,23)(179,23)(205,23)
\put(129,34){\circle*{4}}
\qbezier(154,23)(141,29)(129,34)
\put(129,55){\circle*{4}}
\qbezier(154,65)(141,60)(129,55)
\put(243,97){\circle*{4}}
\qbezier(176,97)(209,97)(243,97)
\put(243,122){\circle*{4}}
\qbezier(243,97)(243,110)(243,122)
\put(215,97){\circle*{4}}
\put(215,121){\circle*{4}}
\qbezier(215,97)(215,109)(215,121)
\put(232,66){\circle*{4}}
\qbezier(205,66)(218,66)(232,66)
\put(232,23){\circle*{4}}
\qbezier(205,23)(218,23)(232,23)
\put(140,70){$v_{1}$}
\put(143,12){$v_{2}$}
\put(197,12){$v_{3}$}
\put(187,62){$v_{4}$}
\put(169,104){$v_{5}$}
\put(166,-17){$G_{16}$}
\put(200,86){$v_{s+5}$}
\put(565,97){\circle*{4}}
\put(546,66){\circle*{4}}
\qbezier(565,97)(555,82)(546,66)
\put(600,68){\circle*{4}}
\qbezier(565,97)(582,83)(600,68)
\put(546,24){\circle*{4}}
\qbezier(546,66)(546,45)(546,24)
\put(600,24){\circle*{4}}
\qbezier(600,68)(600,46)(600,24)
\qbezier(546,24)(573,24)(600,24)
\put(527,80){\circle*{4}}
\qbezier(546,66)(536,73)(527,80)
\put(526,44){\circle*{4}}
\qbezier(546,24)(536,34)(526,44)
\put(628,36){\circle*{4}}
\qbezier(600,24)(614,30)(628,36)
\put(628,75){\circle*{4}}
\qbezier(600,68)(614,72)(628,75)
\put(537,97){\circle*{4}}
\qbezier(565,97)(551,97)(537,97)
\put(633,97){\circle*{4}}
\qbezier(565,97)(599,97)(633,97)
\put(633,122){\circle*{4}}
\qbezier(633,97)(633,110)(633,122)
\put(605,97){\circle*{4}}
\put(605,121){\circle*{4}}
\qbezier(605,97)(605,109)(605,121)
\put(550,59){$v_{1}$}
\put(543,12){$v_{2}$}
\put(591,12){$v_{3}$}
\put(583,61){$v_{4}$}
\put(558,104){$v_{5}$}
\put(561,-17){$G_{19}$}
\put(251,-52){Fig. 4.1. $G_{1}-G_{19}$}
\put(304,97){\circle*{4}}
\put(282,63){\circle*{4}}
\qbezier(304,97)(293,80)(282,63)
\put(336,65){\circle*{4}}
\qbezier(304,97)(320,81)(336,65)
\put(282,22){\circle*{4}}
\qbezier(282,63)(282,43)(282,22)
\put(336,22){\circle*{4}}
\qbezier(336,65)(336,44)(336,22)
\qbezier(282,22)(309,22)(336,22)
\put(367,97){\circle*{4}}
\qbezier(304,97)(335,97)(367,97)
\put(340,97){\circle*{4}}
\put(340,119){\circle*{4}}
\qbezier(340,97)(340,108)(340,119)
\put(367,119){\circle*{4}}
\qbezier(367,97)(367,108)(367,119)
\put(276,97){\circle*{4}}
\qbezier(304,97)(290,97)(276,97)
\put(258,63){\circle*{4}}
\qbezier(282,63)(270,63)(258,63)
\put(257,22){\circle*{4}}
\qbezier(282,22)(269,22)(257,22)
\put(362,22){\circle*{4}}
\qbezier(336,22)(349,22)(362,22)
\put(286,58){$v_{1}$}
\put(277,10){$v_{2}$}
\put(333,10){$v_{3}$}
\put(341,61){$v_{4}$}
\put(300,102){$v_{5}$}
\put(298,-17){$G_{17}$}
\put(435,97){\circle*{4}}
\put(410,66){\circle*{4}}
\qbezier(435,97)(422,82)(410,66)
\put(410,22){\circle*{4}}
\qbezier(410,66)(410,44)(410,22)
\put(466,66){\circle*{4}}
\qbezier(435,97)(450,82)(466,66)
\put(466,22){\circle*{4}}
\qbezier(410,22)(438,22)(466,22)
\qbezier(466,66)(466,44)(466,22)
\put(385,22){\circle*{4}}
\qbezier(410,22)(397,22)(385,22)
\put(386,66){\circle*{4}}
\qbezier(410,66)(398,66)(386,66)
\put(408,97){\circle*{4}}
\qbezier(435,97)(421,97)(408,97)
\put(506,97){\circle*{4}}
\qbezier(435,97)(470,97)(506,97)
\put(473,97){\circle*{4}}
\put(473,119){\circle*{4}}
\qbezier(473,97)(473,108)(473,119)
\put(506,119){\circle*{4}}
\qbezier(506,97)(506,108)(506,119)
\put(497,66){\circle*{4}}
\qbezier(466,66)(481,66)(497,66)
\put(414,59){$v_{1}$}
\put(406,10){$v_{2}$}
\put(469,14){$v_{3}$}
\put(449,59){$v_{4}$}
\put(430,102){$v_{5}$}
\put(425,-17){$G_{18}$}
\end{picture}
\end{center}

\

{\bf Case 2} There is only $1$ $p$-dominator on $\mathbb{C}$ (see $G_{2}-G_{4}$ in Fig. 4.1).

{\bf Subcase 2.1} For $G_{2}$, let $\mathbb{H}=G_{2}-v_{3}v_{4}+v_{3}v_{1}$. As Case 1, it is proved that $\gamma(G_{2})\leq \gamma(\mathbb{H})$ and $q_{min}(\mathbb{H})< q_{min}(G_{2})$.

{\bf Subcase 2.2} For $G_{3}$, suppose $X=(\,x_1$, $x_2$, $\ldots$, $x_{n-1}$, $x_{n}\,)^T$ is a unit eigenvector corresponding to $q_{min}(G_{3})$.

{\bf Claim} $|x_{4}|>|x_{1}|$, $|x_{5}|>|x_{3}|$. Denote by $v_{k}$ the pendant vertex attached to $v_{4}$. Suppose $0<|x_{4}|\leq|x_{1}|$. Let $G^{'}_{3}=G_{3}-v_{4}v_{k}+v_{1}v_{k}$. By Lemma \ref{le0,04}, then $q_{min}(G^{'}_{3})<q_{min}(G_{3})$. This is a contradiction because $G^{'}_{3}\cong G_{3}$. Suppose $|x_{4}|=|x_{1}|=0$. By Lemma \ref{le04,01}, we get that $x_{2}\neq 0$, $x_{3}\neq 0$. By $q_{min}(G_{3})x_{2}=2x_{2}+x_{3}$, $q_{min}(G_{3})x_{3}=2x_{3}+x_{2}$, we get $x^{2}_{2}=x^{2}_{3}$. Suppose $x_{2}>0$. Then we get $q_{min}(G_{3})x_{2}=2x_{2}+x_{3}\geq x_{2}$. This means that $q_{min}(G_{3})\geq 1$ which contradicts $q_{min}(G_{3})<1$. Thus, $|x_{4}|>|x_{1}|$. Similarly, we get $|x_{5}|>|x_{3}|$. Then the claim holds.

Suppose $|x_{1}|=\min\{|x_{1}|$, $|x_{2}|$, $|x_{3}|\}$ and $x_{1}\geq 0$. If $|x_{2}|> |x_{5}|$, by Lemma \ref{le04,01}, suppose $x_1x_5\geq 0$. Let $H=G_{3}-v_{1}v_{5}$. Also by Lemma \ref{le04,01}, suppose for any $j\neq 1, 5$, $\mathrm{sgn}x_j= (-1)^{d_{istH}(v_{j}, v_1)}$. Let  $\mathbb{H}=G_{3}-v_{1}v_{5}+v_{3}v_{1}$. Because $|x_{5}|>|x_{3}|$, it follows that $q_{min}(\mathbb{H})\leq X^{T}Q(\mathbb{H})X<X^{T}Q(G_{3})X =q_{min}(G_{3})$. Let $B_{1}=\mathbb{H}[v_{1}, v_{2}]$, $B_{2}=\mathbb{H}- \{v_{1}, v_{2}\}$. As Lemma \ref{le03,07}, we can get a minimal dominating set $D$ of $\mathbb{H}$, which contains all $p$-dominators but no any pendant vertex and no $v_{3}$, such that $D=\{v_{1}\}\cup D_{2}$, where $D_{2}$ is a dominating set of $B_{2}$. Note that $D$ is also a dominating set of $G_{3}$. So, $\gamma(G_{3})\leq \gamma(\mathbb{H})$. If $|x_{2}|< |x_{5}|$, by Lemma \ref{le04,01}, $x_{1}x_{2}\geq 0$. Let $H=G_{3}-v_{1}v_{2}$. Also by Lemma \ref{le04,01}, suppose for any $j\neq 1, 2$, $\mathrm{sgn}x_j= (-1)^{d_{istH}(v_{j}, v_1)}$. Let  $\mathbb{H}=G_{3}-v_{1}v_{5}+v_{3}v_{1}$. Because $|x_{5}|>|x_{3}|$, it follows that $q_{min}(\mathbb{H})< q_{min}(G_{3})$ similarly. As the case that $|x_{2}|> |x_{5}|$, it is proved that $\gamma(G_{3})\leq \gamma(\mathbb{H})$. If $|x_{2}|= |x_{5}|$, by Lemma \ref{le04,01}, without loss of generality, suppose $x_{1}x_{5}\geq 0$. Let  $\mathbb{H}=G_{3}-v_{1}v_{5}+v_{3}v_{1}$. As the case that $|x_{2}|> |x_{5}|$, it is proved that $q_{min}(\mathbb{H})< q_{min}(G_{3})$, $\gamma(G_{3})\leq \gamma(\mathbb{H})$.

For the both cases that $|x_{2}|=\min\{|x_{1}|$, $|x_{2}|$, $|x_{3}|\}$ and $|x_{3}|=\min\{|x_{1}|$, $|x_{2}|$, $|x_{3}|\}$. As the case that $|x_{1}|=\min\{|x_{1}|$, $|x_{2}|$, $|x_{3}|\}$, it is proved that there exists a graph $\mathbb{H}$ such that $g(\mathbb{H})=3$, $\gamma(G_{3})\leq \gamma(\mathbb{H})$ and $q_{min}(\mathbb{H})< q_{min}(G_{3})$.

In a same way, for $G_{4}$, it is proved that there exists a graph $\mathbb{H}$ such that $g(\mathbb{H})=3$, $\gamma(G_{4})\leq \gamma(\mathbb{H})$ and $q_{min}(\mathbb{H})< q_{min}(G_{4})$.

And in a same way, for the cases that {\bf Case 3} there is exactly $2$ $p$-dominators on $\mathbb{C}$ (see $G_{5}-G_{10}$ in Fig. 4.1); {\bf Case 4} there is exactly $3$ $p$-dominators on $\mathbb{C}$ (see $G_{11}-G_{15}$ in Fig. 4.1); {\bf Case 5} there is exactly $4$ $p$-dominators on $\mathbb{C}$ (see $G_{16}-G_{18}$ in Fig. 4.1); {\bf Case 6} there is exactly $5$ $p$-dominators on $\mathbb{C}$ (see $G_{19}$ in Fig. 4.1), it is proved that the exists a a graph $\mathbb{H}$ such that $g(\mathbb{H})=3$, $\gamma(G)\leq \gamma(\mathbb{H})$ and $q_{min}(\mathbb{H})\leq q_{min}(G)$.
Thus, the result follows as desired. \ \ \ \ \ $\Box$
\end{proof}

\setlength{\unitlength}{0.6pt}
\begin{center}
\begin{picture}(355,114)
\put(3,96){\circle*{4}}
\put(3,19){\circle*{4}}
\qbezier(3,96)(3,58)(3,19)
\put(49,58){\circle*{4}}
\qbezier(3,96)(26,77)(49,58)
\qbezier(3,19)(26,39)(49,58)
\put(187,58){\circle*{4}}
\put(111,58){\circle*{4}}
\put(209,58){\circle*{4}}
\put(198,58){\circle*{4}}
\put(226,58){\circle*{4}}
\put(226,88){\circle*{4}}
\put(288,58){\circle*{4}}
\qbezier(226,58)(226,73)(226,88)
\put(170,58){\circle*{4}}
\put(347,58){\circle*{4}}
\put(1,102){$v_{2}$}
\put(47,45){$v_{3}$}
\put(102,45){$v_{4}$}
\put(1,6){$v_{1}$}
\put(353,83){$v_{\frac{n+5}{2}}$}
\put(335,43){$v_{\frac{n+3}{2}}$}
\qbezier(3,19)(3,58)(3,96)
\qbezier(3,19)(3,58)(3,96)
\qbezier(3,19)(3,58)(3,96)
\qbezier(3,19)(3,58)(3,96)
\put(347,85){\circle*{4}}
\qbezier(347,58)(347,72)(347,85)
\qbezier(288,58)(317,58)(347,58)
\qbezier(288,58)(257,58)(226,58)
\put(288,88){\circle*{4}}
\qbezier(288,58)(288,73)(288,88)
\put(111,87){\circle*{4}}
\qbezier(111,58)(111,73)(111,87)
\qbezier(111,58)(140,58)(170,58)
\qbezier(111,58)(80,58)(49,58)
\put(170,87){\circle*{4}}
\qbezier(170,58)(170,73)(170,87)
\put(277,45){$v_{\frac{n+1}{2}}$}
\put(129,-9){Fig. 4.2. $\mathscr{H}_{3,\frac{n-3}{2}}$}
\end{picture}
\end{center}

\begin{lemma}\label{le04,06} 
Let $G$ be a nonbipartite $\mathcal {F}_{g, l}$-graph of order $n$ for some $l$ and with domination number $\frac{n-1}{2}$. Then $q_{min}(G)\geq q_{min}(\mathscr{H}_{3,\frac{n-3}{2}})$ with equality if and only if $G\cong \mathscr{H}_{3,\frac{n-3}{2}}$ (see Fig. 4.2).
\end{lemma}

\begin{proof}
Because $G$ is nonbipartite, $g$ is odd. If $G$ is a $\mathcal {F}^{\circ}_{g, l}$-graph, then the theorem follows from Lemma \ref{th04,05}. If $g=3$, then the theorem follows from Theorem \ref{th04,04}. For $g=5$, the theorem follows from Lemma \ref{le04,09}. Next we consider the case that $G$ is not a $\mathcal {F}^{\circ}_{g, l}$-graph and suppose $g\geq 7$.

Let $X=(\,x_{1}$, $x_{2}$, $\ldots$, $x_{n}\,)^T$ is a unit eigenvector corresponding to $q_{min}(G)$. Suppose $x_{a}=\min\{|x_{1}|$, $|x_{2}|$, $\ldots$, $|x_{g}|\}$. Note that by Theorem \ref{le03,08}, in $G$, there are at most $3$ consecutive vertices of $\mathbb{C}$ such that none of them is $p$-dominator, and there are $2$ cases as follows to consider.

{\bf Case 1} In $G$, there is exactly one vertex of $\mathbb{C}$ which is not $p$-dominator. Note that $G$ is not a $\mathcal {F}^{\circ}_{g, l}$-graph. Then $n\geq g+2$ and $v_{g}$ is the only one vertex which is not $p$-dominator on $\mathbb{C}$. By a same discussion in the proof of Lemma \ref{le04,03} (see \cite{YMCG}), it is proved that $x_{g}=\max\{|x_{1}|$, $|x_{2}|$, $\ldots$, $|x_{g-1}|$, $|x_{g}|\}$. Then we suppose $a\leq g-1$. By Lemma \ref{le04,01}, if $a\leq g-3$, without loss of generality, suppose $x_{a+1}\leq x_{a-1}$, $x_{a+1}x_{a}\geq 0$, $|x_{a-1}|\geq |x_{a+2}|$. Let $G_{1}=G-v_{a}v_{a-1}+v_{a}v_{a+2}$ (if $|x_{a-1}|\leq |x_{a+2}|$ and $a\geq 2$, let $G_{1}=G-v_{a+1}v_{a+2}+v_{a+1}v_{a-1}$; if $a=1$, let $G_{1}=G-v_{1}v_{g}+v_{1}v_{3}$). If $a= g-2$, suppose $|x_{g-1}|\leq |x_{g-3}|$, $x_{g-1}x_{g-2}\geq 0$, and then let $G_{1}=G-v_{g-1}v_{g}+v_{g-1}v_{g-3}$. If $a= g-1$, because $|x_{g}|\geq |x_{g-2}|$, then suppose $x_{g-1}x_{g-2}\geq 0$. Let $G_{1}=G-v_{g-1}v_{g}+v_{g-1}v_{g-3}$. Note that $\gamma(G_{1})\leq \frac{n-1}{2}$. As the proof of Lemma \ref{th04,02}, we get that $\gamma(G)\leq \gamma(G_{1})= \frac{n-1}{2}$, $q_{min}(G_{1})< q_{min}(G)$. Note that $g(G_{1})=3$. Then the theorem follows from Theorem \ref{th04,04}.

{\bf Case 2} In $G$, there are exactly $3$ consecutive vertices of $\mathbb{C}$ such that each of them is not $p$-dominator.
Note that $G$ is not a $\mathcal {F}^{\circ}_{g, l}$-graph. Combined with Theorem \ref{le03,08}, the $3$ vertices of $\mathbb{C}$ such that each of them is not $p$-dominator are $v_{g-2}$, $v_{g-1}$, $v_{g}$ or $v_{g}$, $v_{1}$, $v_{2}$. Without loss of generality, we suppose the $3$ vertices are $v_{g-2}$, $v_{g-1}$, $v_{g}$. By Lemma \ref{th03,02}, $|x_{g}|>0$. We say that $|x_{g}|>|x_{g-2}|$. Otherwise, suppose $|x_{g}|\leq|x_{g-2}|$. Let $G^{'}=G-v_{g}v_{g+1}+v_{g+1}v_{g-2}$. Then by Lemma \ref{le0,04}, $q_{min}(G^{'})< q_{min}(G)$. This is a contradiction because $G^{'}\cong G$. Hence $|x_{g}|>|x_{g-2}|$. And then $a\leq g-1$.

{\bf Subcase 2.1} $a\leq g-4$. By Lemma \ref{le04,01}, without loss of generality, suppose $x_{a+1}\leq x_{a-1}$, $x_{a+1}x_{a}\geq 0$.
As Case 1, it is proved that the theorem holds.

{\bf Subcase 2.2} $a= g-3$. By Lemma \ref{le04,01}, suppose $x_{g-2}\leq x_{g-4}$, $x_{g-2}x_{g-3}\geq 0$; suppose $|x_{g-4}|\geq |x_{g-1}|$. Denote by $v_{\tau_{g-3}}$ the pendant vertex attached to $v_{g-3}$. Let $G_{1}=G-v_{g-3}v_{g-4}+v_{g-3}v_{g-1}-v_{g-3}v_{\tau_{g-3}}+v_{g}v_{\tau_{g-3}}$ (if $x_{g-4}\leq x_{g-1}$, let $G_{1}=G-v_{g-2}v_{g-1}+x_{g-2}x_{g-4}$). As Case 1, it is proved that the theorem holds.

{\bf Subcase 2.3} $a= g-2$. By Lemma \ref{le04,01}, suppose $x_{g-1}\leq x_{g-3}$, $x_{g-1}x_{g-2}\geq 0$; suppose $|x_{g-3}|\geq |x_{g}|$. Denote by $v_{\tau_{g-3}}$ the pendant vertex attached to $v_{g-3}$. Let $G_{1}=G-v_{g-2}v_{g-3}+v_{g-2}v_{g}$ (if $x_{g-3}\leq x_{g}$, let $G_{1}=G-v_{g-1}v_{g}+x_{g-1}x_{g-3}-v_{g-3}v_{\tau_{g-3}}+v_{g}v_{\tau_{g-3}}$). As Case 1, it is proved that the theorem holds.

{\bf Subcase 2.4} $a= g-1$. Note $|x_{g}|>|x_{g-2}|$. By Lemma \ref{le04,01}, $x_{g-2}x_{g-1}\geq 0$. Without loss of generality, suppose $x_{g-3}\geq x_{g}$, let $G_{1}=G-v_{g-2}v_{g-3}+v_{g-2}v_{g}$ (if $x_{g-3}\leq x_{g}$, let $G_{1}=G-v_{g-1}v_{g}+x_{g-1}x_{g-3}-v_{g}v_{g+1}+v_{g-3}v_{g+1}$). As Case 1, it is proved that the theorem holds.
This completes the proof. \ \ \ \ \ $\Box$
\end{proof}

By Lemmas \ref{th03,02}, \ref{le04,06}, we get the following Theorem \ref{th04,07}.

\begin{theorem}\label{th04,07} 
Let $G$ be a nonbipartite connected unicyclic graph of order $n\geq 3$ and with domination number $\frac{n-1}{2}$. Then $q_{min}(G)\geq q_{min}(\mathscr{H}_{3,\frac{n-3}{2}})$ with equality if and only if $G\cong \mathscr{H}_{3,\frac{n-3}{2}}$.
\end{theorem}

\section{Proof of main results}

\begin{Proof}
By Lemmas \ref{le02,01}, \ref{le0,07}, then $G$ contains a nonbipartite
unicyclic spanning subgraph $H$
with $g_{o}(H)=g_{o}(G)$, $\gamma(H)=\gamma(G)$ and $q_{min}(H)\leq q_{min}(G)$. By Theorem \ref{th04,07}, it follows that $q_{min}(H)\geq q_{min}(\mathscr{H}_{3,\frac{n-3}{2}})$ with equality if and only if $H\cong \mathscr{H}_{3,\frac{n-3}{2}}$. Thus it follows that $q_{min}(G)\geq q_{min}(\mathscr{H}_{3,\frac{n-3}{2}})$.

Suppose that $q_{min}(G)= q_{min}(\mathscr{H}_{3,\frac{n-3}{2}})$. Then $q_{min}(H)= q_{min}(\mathscr{H}_{3,\frac{n-3}{2}})$ and $H\cong \mathscr{H}_{3,\frac{n-3}{2}}$. For convenience, we suppose that $H= \mathscr{H}_{3,\frac{n-3}{2}}$.
Suppose that $Y$ is a unit eigenvector corresponding to $q_{min}(G)$. Note that $q_{min}(\mathscr{H}_{3,\frac{n-3}{2}})=q_{min}(H)\leq Y^{T}Q(H)Y\leq Y^{T}Q(G)Y=q_{min}(G)$. Because we suppose that $q_{min}(G)= q_{min}(\mathscr{H}_{3,\frac{n-3}{2}})$, it follows that $Y^{T}Q(H)Y= Y^{T}Q(G)Y$ and $Q(H)Y=q_{min}(H)Y$.

For $\mathscr{H}_{3,\frac{n-3}{2}}$ (see Fig. 4.2), we claim that $y_{3}>y_{1}$, $y_{3}>y_{2}$. Otherwise, suppose that $y_{3}\leq y_{1}$. Let $H'=\mathscr{H}_{3,\frac{n-3}{2}}-v_{3}v_{4}+v_{1}v_{4}$. By Lemma \ref{le0,04}, it follows that $q_{min}(H')< q_{min}(\mathscr{H}_{3,\frac{n-3}{2}})$. This is a contradiction because $H'\cong H\cong \mathscr{H}_{3,\frac{n-3}{2}}$. Thus our claim holds.

If $G\neq H$, combined with Lemma \ref{le02,03}, then for any edge $v_{i}v_{j}\not\in E(H)$, it follows that $x_{i}+x_{j}\neq 0$, and then $Y^{T}Q(H)Y< Y^{T}Q(G)Y$, which contradicts $Y^{T}Q(H)Y= Y^{T}Q(G)Y$. Then it follows that $q_{min}(G)= q_{min}(\mathscr{H}_{3,\frac{n-1}{2}})$ if and only if $G\cong \mathscr{H}_{3,\frac{n-1}{2}}$.
This completes the proof. \ \ \ \ \ $\Box$
\end{Proof}

In a same way, with Lemmas \ref{th03,05},  \ref{cl03,06} and \ref{le04,09}, Theorem \ref{th05,4} is proved.

{\bf Remark} It can be seen that the conjecture in \cite{YMCG} that $\mathbb{S}$ has the smallest $q_{min}$ holds for the graphs with domination number $\gamma =\frac{n-1}{2}$ and the graphs with girth at most 5. With references \cite{YGZW} and
\cite{YMCG}, it can also be seen that the minimum $q_{min}$ of the connected nonbipartite graph on $n\geq 5$ vertices, with domination number $\frac{n+1}{3}<\gamma\leq \frac{n-2}{2}$ and girth $g\geq 5$, is still open.

\small {

}

\end{document}